\documentclass[reqno]{amsart}

\usepackage[utf8]{inputenc}
\usepackage[english]{babel}
\usepackage[margin=1.6in]{geometry}
\usepackage[T1]{fontenc}
\usepackage{amsmath}
\theoremstyle{plain}
\usepackage{amssymb}
\usepackage{amsthm}
\usepackage{graphics}
\usepackage{euler}
\usepackage{calrsfs}
\usepackage{palatino}
\usepackage{graphicx}
\usepackage{amsmath}
\usepackage{amsfonts}
\usepackage{amssymb}
\usepackage{amsthm}
\usepackage{amstext}
\usepackage{amscd}
\usepackage{mathtools} 
\usepackage{hyperref}
\usepackage{tikz}
\usepackage{tikz-cd}
\usetikzlibrary{matrix}
\usepackage[all]{xy}
\usepackage{algorithm}
\usepackage{algorithmic}
\usetikzlibrary{patterns,decorations.pathreplacing,babel}
\usetikzlibrary{shapes,intersections}

\usepackage{hyperref}
\hypersetup{
    colorlinks=true,
    linkcolor=purple,
    filecolor=magenta,      
    citecolor=cyan,
    }

\newcommand\C{\mathbb{C}}

\renewcommand\P{\mathbb{P}}

\newcommand{\I}{\mathcal{I}}

\newcommand{\OO}{\mathcal{O}}
\renewcommand{\I}{\mathcal{I}}
\renewcommand{\L}{\mathcal{L}}

\newcommand{\U}{\mathcal{U}}

\newcommand{\X}{\mathcal{X}}

\newcommand{\F}{\mathscr{F}}
\newcommand{\PP}{\mathcal{P}}
\newcommand{\E}{\mathcal{E}}

\newcommand{\dF}{\mathcal{F}}

\newcommand{\D}{\mathscr{D}}

\renewcommand\v{\raise0.9ex\hbox{$\scriptscriptstyle\vee$}}

\DeclareMathOperator{\Spec}{Spec}
\DeclareMathOperator{\Ex}{Ex}
\DeclareMathOperator{\Der}{Der}
\DeclareMathOperator{\Ext}{Ext}

\DeclareMathOperator{\Sing}{Sing}
\DeclareMathOperator{\codim}{codim}

\DeclareMathOperator{\Per}{Per}

\DeclareMathOperator{\Lef}{Lef}
\DeclareMathOperator{\Leff}{Leff}
\DeclareMathOperator{\Int}{int}

\DeclareFontFamily{OT1}{pzc}{}
\DeclareFontShape{OT1}{pzc}{m}{it}{<-> s * [1.100] pzcmi7t}{}
\DeclareMathAlphabet{\mathpzc}{OT1}{pzc}{m}{it}

\usepackage{hyperref}

\newcommand{\mybar}[3]{%
    \mathrlap{\hspace{#2}\overline{\scalebox{#1}[1]{\phantom{\ensuremath{#3}}}}}\ensuremath{#3}
}

\newtheorem{theorem}{Theorem}[section]
\newtheorem{proposition}[theorem]{Proposition}
\newtheorem{lemma}[theorem]{Lemma}
\newtheorem{corollary}[theorem]{Corollary}

\theoremstyle{definition}

\newtheorem{definition}[theorem]{Definition}
 \newtheorem{remark}[theorem]{Remark}

\newtheorem{example}[theorem]{Example}

\newtheoremstyle{named}{}{}{\itshape}{}{\bfseries}{.}{.5em}{\thmnote{#3}}

\theoremstyle{named}

\title{Extendability of foliations}
\author{Pablo Perrella$^1$ and Sebastian Velazquez$^2$}

\thanks{$^1$ The author was supported by CONICET, Argentina.}
\thanks{$^2$The author was supported by EPSRC, United Kingdom.}
\date{}

\address{Departamento de Matem\'atica 
, Facultad de Ciencias Exactas y Naturales, Universidad de Buenos Aires 
and IMAS (UBA-CONICET) 
, Buenos Aires, Argentina.}
\email{pperrella@dm.uba.ar}

\address{Department of Mathematics, King's College London, Strand,
London WC2R 2LS, United Kingdom.}
\email{sebastian.velazquez@kcl.ac.uk}

\keywords{Holomorphic foliations, unfoldings, extension of foliations, tubular neighborhoods}
\subjclass[2020]{14D15, 14D15, 14D06, 	32S65, 32M25}

\begin{document}

\begin{abstract}
       Given a foliation $\F$ on $X$ and an embedding $X\subseteq Y$, is there a foliation on $Y$ extending $\F$? Using formal methods, we show that this question has an affirmative answer whenever the embedding is sufficiently positive with respect to $(X,\F)$ and the singularities of $\F$ belong to a certain class. These tools also apply in the case where $Y$ is the total space of a deformation of $X$. Regarding the uniqueness of the extension, we prove a foliated version of a statement by Fujita and Grauert ensuring the existence of tubular neighborhoods. We also give sufficient conditions for a foliation to have only trivial unfoldings, generalizing a result due to G\'omez-Mont.
\end{abstract}

\maketitle

\section{Introduction}
In the study of any geometric structure, a good understanding of the process of restriction/extension with respect to a subvariety $X\subseteq Y$ provides a great deal of clarity to the problem. This article is dedicated to the study of the extension problem for holomorphic singular foliations in both the analytic and algebraic settings. Broadly speaking, this involves determining whether, given a foliation $\F$ on $X$ and an embedding $X\subseteq Y$, there exists a foliation on $Y$ that restricts to $\F$ on $X$. In the case where $Y=\P^m$, this question was already raised in \cite[Section 2.1]{cerveau2013quelques} and addressed in recent works like \cite{araujo2018codimension,   figueira2023extensions, benedetti2023codimension}. In this work we take a different approach that allows us to work on more general settings. 

The motivation for this work has two facets which, although complementary, deserve to be mentioned separately. On the one hand, the study of the geometry of the restriction of a foliation to a suitable subvariety has been a key point throughout the literature. Additionally, it is noteworthy that, with the appropriate formulation, a special case of this issue is the theory of unfoldings. An unfolding of a foliation $\F$ on a variety $X$ over some base $T$ is a deformation $\X\to T$ of $X$, together with a foliation on $\X$ restricting to $\F$ on the central fiber. 
This framework was introduced by T. Suwa in \cite{SuwaVersality} and has remained an active area of research ever since (see for instance \cite{gomez1989unfoldings, molinuevo2016unfoldings, persistent2021}. The theory of unfoldings is closely related to the deformation theory of foliations, since of course the fibers of an unfolding define a family of foliations. 
In particular, the content of this work can also be intended as a contribution towards a better comprehension of the moduli theory of foliations.
Let us now briefly explain our approach.

 Of course, if we expect to be able to construct an extension, the singularities of $\F$ should be, in some way, locally unobstructed. Given a germ of foliation $\F$, an unfolding of order $n$ of $\F$ is an unfolding with base $\Spec(\C[t]/(t^{n+1}))$. We say that a germ $\F$ is  \emph{unobstructed} if every unfolding of order $n$ can be extended to higher order. We say that $\F$ is \emph{rigid} if every first order unfolding is isomorphic to the trivial unfolding. 
A special case of foliation singularities are the so-called \emph{persistent} singularities defined in \cite{persistent2021}. These are the singularities $p\in \Sing(\F)$ that cannot be smoothed out 
along an unfolding (see Definition \ref{persistent} below). We will denote by $\I_{\Per}$ the sheaf of ideals associated to this type of singularities.

The starting point of our strategy is to understand the extension problem at the level of every infinitesimal neighborhood of $X$ in $Y$, describing the groups where the obstructions and the different extensions of $\F$ naturally reside. As a consequence, 
we will be able to show that $\F$ admits a \emph{formal} extension to $Y$ whenever $X\subseteq Y$ is a codimension $1$ regular embedding (analytical or algebraic), $\F$ has unobstructed singularities and 
$H^1\big(X,\I_{\Per}\otimes N_\F\otimes N_{X/Y}^{\otimes -n}\big)=0$
for every $n\geqslant 1$, where $N_\F$ is the normal bundle of the foliation. Once the formal problem is settled, we will rely on the effective Lefschetz condition $\Leff(X,Y)$ (Definition \ref{lefschetz}) in order to give a positive answer to the extension problem:

\begin{theorem}
   Let $\F$ be a codimension $1$ foliation on a smooth variety $X$ having unobstructed singularities and let $X\subseteq Y$ be a regular embedding with $\dim(Y)=\dim(X)+1$.  
   Suppose that $$H^1(X,\I_{\Per}\otimes N_\F\otimes N_{X/Y}^{\otimes -n})=0$$
   for every $n\geqslant 1$ and $\Leff(X,Y)$ holds.
    Then $\F$ extends to a foliation on $Y$. 
\end{theorem}

This result can be applied to the following two cases: (i) $X\subseteq Y$ is an ample divisor; and (ii) $X\subseteq \X$ is the central fiber of a deformation. Hence, the previous theorems can be further refined to obtain the following conclusions.

\begin{corollary}
     Let $X$ be a smooth ample divisor on a smooth projective variety $Y$ and $\F$ a codimension $1$ foliation on $X$ with unobstructed singularities and $H^1\big(X,\I_{\Per}\otimes N_\F\otimes N_{X/Y}^{\otimes - n}\big)=0$ for every $n\geqslant 1$. Then $\F$ extends to a foliation on $Y$.
\end{corollary}

\begin{corollary}
    Let $\F$ be a codimension $1$ foliation on a proper smooth variety $X$ having unobstructed singularities and such that $H^1(X,\I_{\Per}\otimes N_\F)=0.$ Then $\F$  unfolds along any deformation $\X\to \Spec(\C[[t]])$ of $X$.
\end{corollary}

If, on the contrary, one seeks to extend a foliation under mild assumptions on its singularities, one should rely on a stronger positivity condition for the embedding. This is the strategy applied in the previous works mentioned earlier in the introduction. In \cite{figueira2023extensions}, for instance, the author shows that a foliation $\F$ on a hypersurface $X\subseteq \P^m$ extends whenever $\deg(X)> 2 \deg(\F) +1$. In Example \ref{ej:figueira} we recover this statement for foliations with unobstructed singularities. This is a special case of the following: 

\begin{theorem}
Let $X$ be a smooth ample divisor on a smooth projective variety $Y$ and $\F$ a codimension $1$ foliation on $X$ with unobstructed singularities. If $\dim X\geqslant 3$ and $N_\F^{\otimes 2}\otimes N_{X/Y}^{\otimes -1}$ is anti-ample, then $\F$ extends to a unique foliation on $Y$. 
\end{theorem}

One reason a foliation might extend is the existence of a tubular neighborhood, i.e., a neighborhood $U\subseteq Y$ of $X$ endowed with a morphism $U\to X$ restricting to the identity on $X$. The construction of such a neighborhood was achieved independently by Grauert and Fujita in the analytic and algebraic settings respectively, under the common vanishing hypotheses of the groups
\[
H^1\big(X,T_X\otimes N_{X/Y}^{\otimes -n}\big) = 0
\]
for all $n\geqslant 1$ and the ampleness of $X\subseteq Y$ \cite{grauert1962modifikationen, fujita1983rational}. In the same spirit, L'vovsky proved that given a smooth subvariety $X\subseteq \P^m$ such that $H^1(X,T_X(-1)) = 0$, an extension of $X$ along a linear embedding $\P^m\subseteq \P^{m+1}$ must be a cone \cite{l1992extensions}. In this direction, the existence of a special singular foliation on $X$ can be used to ensure the existence of a tubular neighborhood.

\begin{theorem}
  Let $X$ be a smooth projective variety with $\dim(X)\geq 3$ and $X\subseteq Y$ a regular embedding with ample normal bundle. Suppose further that there exists a codimension one foliation $\F$ on $X$ with rigid singularities admitting an extension to $Y$ and satisfying
    $$H^1\big(X,T_\F\otimes N_{X/Y}^{\otimes -n}\big) = 0$$
for each $n\geqslant 1$. Then $X$ has a tubular neighborhood. If moreover $Y$ is projective and $X$ is an ample divisor, then $X$ admits a Zariski tubular neighborhood and we have the following dichotomy:
\begin{enumerate}
    \item the extension of $\F$ is unique, or
    \item $X\subseteq Y$ is a hyperplane included in a projective space and the extension of $\F$ is a cone.
\end{enumerate}
\end{theorem}

This result admits an analogous version in the unfolding scenario. One of the earliest versions can be found in \cite{gomez1989unfoldings}. It was proved that given a compact surface $X$ and a foliation $\F$ by curves on $X$ such that $H^1(X,T_\F)= 0$ and whose singular points have multiplicity one and non-vanishing trace, the only unfoldings of $\F$ are (analytically) the trivial ones. The second manifestation appears on \cite{massri2018kupka}, where the same conclusion was achieved for codimension one foliations $\F$ on $\P^m$ whose tangent sheaf splits as a sum of line bundles and whose non-Kupka locus has codimension at least three. As we will see later (Lemma \ref{ex:locallyfree}), these singularities are rigid and hence the following is a common generalization of both results.

\begin{theorem}
Let $\F$ be a codimension $1$ foliation with rigid singularities on a proper smooth variety $X$. If $H^1(X,T_\F) = 0$, then every unfolding of $\F$ is analytically trivial. 
\end{theorem}
\,

The paper is organized as follows. In Section \ref{section2}, we gather definitions and preliminary results that will play a central role in our approach.
Section 3 is devoted to the study of the extension problem over infinitesimal neighborhoods.
In \ref{section4.1}, we apply these results to prove Theorem 1.1 and derive Corollaries 1.2 and 1.3. In \ref{section4.2}, we establish a positivity criterion for extendability, which corresponds to the first assertion of Theorem 1.4. The second part of that statement follows from the analysis carried out in Section \ref{section5}, which is devoted to the study of the uniqueness of extensions.
Finally, in Section \ref{section6}, we integrate the techniques developed throughout the paper to study the existence of tubular neighborhoods, culminating in Theorem 1.5 and Theorem 1.6.

\,
\,

\textbf{Acknowledgements.} The work presented in this paper started during a visit of the first author to King's College London, funded by the London Mathematical Society. We are grateful to C. Spicer for making the visit possible and also for many helpful conversations. We would also like to thank the anonymous referee for their careful reading and valuable comments. Gratitude is also due to J.V. Pereira for his comments and suggestions.

\section{Preliminaries}\label{section2}
We will now introduce the notions regarding foliated varieties and their extensions that will be needed for the rest of the article. Throughout this paper we will always assume that $X$ and $Y$ are smooth varieties over the complex numbers, either analytic or algebraic. We choose to work with the terminology of the algebraic setting, but everything carries over verbatim to the analytic category (with schemes replaced by complex spaces, affine open covers by Stein open covers, etc.). Unless stated otherwise, we will not assume them to be compact.

\subsection{Extensions of varieties}

\begin{definition}
 Let $X\subseteq Y$ be a closed embedding. The $n$-th infinitesimal neighborhood $X\subseteq X^{(n)}\subseteq Y$ is the closed subscheme with ideal sheaf $I^{(n)}=I_{X/Y}^{n+1}$.
\end{definition}

For $k<m$ these schemes relate through the closed embeddings $i_{km}:X^{(k)}\hookrightarrow X^{(m)}$ corresponding to the exact sequences
\begin{equation}\label{eq:infneighb}
    0 \rightarrow I_{X/Y}^{k+1}/(I_{X/Y}^{m+1}) \rightarrow \OO_{X^{(m)}}\rightarrow \OO_{X^{(k)}} \rightarrow 0.
\end{equation} 

In particular, we see that for $k=0$ and $m=1$ the leftmost term is a square zero ideal and therefore carries a natural structure of $\OO_X$-module compatible with the quotient map $\OO_{X^{(1)}}\to \OO_X$.

\begin{remark} \label{los_t's}
    We will be particularly interested in the case of codimension $1$ regular embeddings. In this case, the ideal sheaf $I_{X/Y}\subseteq \OO_Y$ is a line bundle, and the exact sequence of structure sheaves for the first infinitesimal neighborhood takes the form
$$ 0 \rightarrow N_{X/Y}^{\otimes -1} \rightarrow \OO_{X^{(1)}}\rightarrow \OO_X \rightarrow 0.$$
If $\U = \{U_i\}$ is a cover of $X$ by affine open subsets such that \smash{$N_{X/Y}^{\otimes -1}\big|_{U_i}\simeq \OO_{U_i}$} together with transition functions $f_{ij}:\OO_{U_i}\vert_{U_{ij}}\to \OO_{U_j}\vert_{U_{ij}}$, we can find local generators $t_i$ of $I_{X/Y}$ such that $t_j=f_{ji} t_i$ as elements of $\OO_{X^{(1)}}$.
More generally, for successive neighborhoods we have exact sequences of the form
\begin{equation}\label{sucesive_thickenings}
0 \rightarrow N_{X/Y}^{\otimes - n} \rightarrow \OO_{X^{(n)}}\rightarrow \OO_{X^{(n-1)}}\rightarrow 0. 
\end{equation}
In terms of the local sections mentioned above, this is just to say that the kernel of these quotient maps is locally generated by the class of the sections $t_i^n$, that glue in $\OO_{X^{(n)}}$ according to the transition functions $\{f_{ji}^n\}$. 
\end{remark}

As mentioned in the introduction, our approach will first take us through a study of the extension problem at the level of the \emph{formal neighborhood} of $X$ in $Y$. We will now explain the corresponding technical aspects needed for our purposes. For a wider exposition on formal neighborhoods the reader is referred to \cite[Section II.9]{hartshorne} or \cite[Section 10.8]{EGAI}.

\begin{definition} Let $X\subseteq Y$ be a closed embedding. The formal neighborhood $\widehat{X}_{/Y}$ of $X$ in $Y$ is the ringed space $(X,\OO_{\widehat{X}_{/Y}}:=\varprojlim \OO_{X^{(n)}})$, where the inverse limit is taken with respect to the homomorphisms appearing in (\ref{eq:infneighb}).
\end{definition}

If $\mathcal{G}$ is a sheaf of $\OO_Y$-modules, then its pullback via the embedding $X^{(n)}\hookrightarrow Y$ defines a sheaf $\mathcal{G}^{(n)}$ on $X^{(n)}$ for every $n$. Moreover, these are equipped with natural restriction maps $r_{km}:\mathcal{G}^{(m)}\to i_{km,*}\mathcal{G}^{(k)}$ (from now on, we will omit the pushforward notation when no confusion arises). Of course, these arrows satisfy $r_{km}=r_{m'm}r_{km'}$ whenever $k<m'<m$. This is a special case of the following:

\begin{definition}
   Let $X\subseteq Y$ be a closed embedding. A (coherent) sheaf $\mathcal G$ on $\widehat{X}_{/Y}$ is an $\OO_{\widehat{X}_{/Y}}$-module of the form $\mathcal G=\varprojlim \mathcal G^{(n)}$, where the $\mathcal G^{(n)}$ is a (coherent) $\OO_{X^{(n)}}$-module and we have for each $k<m$ surjective maps $\varphi_{km}:\mathcal G^{(m)}\to\mathcal G^{(k)}$ of $\OO_{X^{(m)}}$-modules satisfying $\varphi_{km}=\varphi_{m'm}\varphi_{km'}$.
\end{definition}

Understanding the possible extensions of a line bundle to the formal neighborhood of $X$ will be essential for our purposes.
The following statement can be found in \cite[Section I.5.a]{griffiths1966extension}.

\begin{proposition}{\label{ext_line_bundles}}
Let $X^{(n-1)}\subseteq X^{(n)}$ be successive neighborhoods of a codimension $1$ regular embedding $X\subseteq Y$. If $\L^{(n-1)}$ is an extension of a line bundle $\L$ on $X$ to $X^{(n-1)}$, then the vector space
\[
H^1\big(X,N_{X/Y}^{\otimes -n}\big)
\]
acts freely and transitively on the set of extensions of $\L^{(n-1)}$ to $X^{(n)}$.
\end{proposition}
\begin{proof}
Let $\U = \{U_i\}$ be an affine cover of $X$ and $\{f_{ij}\}$ the transition functions of \smash{$N_{X/Y}^{-1}$}. Choose local generators $t_i^n$ of $N_{X/Y}^{\otimes -n}|_{U_i}$ satisfying the equations $t_j^n = f_{ji}^{n}t_i^n$, and transition functions \smash{$\{g_{ij}^{(n-1)}\}$} of $\L^{(n-1)}$. Given two extensions $\L^{(n)}$ and \smash{$\mybar{1}{0pt}{\L}^{(n)}$} of \smash{$\L^{(n-1)}$}, we can find transition functions \smash{$\{g_{ij}^{(n)}\}$} and \smash{$\{\overline{g}_{ij}^{(n)}\}$} that extend \smash{$\{g_{ij}^{(n-1)}\}$} respectively. The exact sequence (\ref{sucesive_thickenings}) implies that
\begin{equation} \label{eq:1}
\overline{g}_{ij}^{(n)} = g_{ij}^{(n)}(1+u_{ij} t_i^n)
\end{equation}
for some local functions $u_{ij}$ on $X$. Combining this with the compatibility relations
\[
g_{ij}^{(n)}g_{jk}^{(n)} = g_{ik}^{(n)}\quad\text{and}\quad \overline{g}_{ij}^{(n)}\overline{g}_{jk}^{(n)} = \overline{g}_{ik}^{(n)}
\] 
we see that
$$(1+u_{ij} t_i^n)(1+u_{jk} t_j^n)=(1+u_{ik} t_i^n).$$
Since $t_j^n=f_{ji}^nt_i^n$, the cocycle conditions
\[
u_{ij}+f_{ji}^{n}u_{jk}-u_{ik} = 0
\]
hold and therefore the collection $\{u_{ij}\}$ defines a class in \smash{$H^1\big(X,N_{X/Y}^{\otimes -n}\big)$} that does not depend on the choices made along the way. The action mentioned earlier is defined by means of the equations (\ref{eq:1}). 
\end{proof}

\subsection{Extensions of distributions and foliations}

\begin{definition}
A codimension $q$ \textit{distribution} $\F$ over a smooth variety $X$ is a coherent subsheaf $I_\F\subseteq \Omega^1_X$
such that $\Omega^1_\F:=\Omega^1_X/I_\F$ is torsion free and the generic rank of $I_\F$ is $q$. The sheaf $\Omega^1_\F$ is called the \textit{conormal} sheaf of the distribution. We say that $\F$ is a \textit{foliation} if through a general point of $X$ there are local generators $\omega_1,\dots,\omega_q$ of $I_\F$ satisfying the \textit{integrability conditions}
\[
\omega_1\wedge\cdots\wedge \omega_q \wedge d\omega_i = 0 \quad\quad (1\leqslant i \leqslant q).
\]
\end{definition}
In the codimension $1$ case, the \textit{normal bundle} $N_\F = I_\F^\vee$ is a line bundle and therefore the distribution $\F$ is uniquely determined by a twisted $1$-form 
$$\omega_\F\in H^0\big(X,\Omega^1_X\otimes N_\F\big)$$ 
up to a global unit of $\OO_X$. The torsion-freeness of $\Omega^1_\F$ is equivalent to $\omega_\F$ not vanishing along a divisor, and the integrability conditions are reduced to the single equation
\[
\omega_\F\wedge d\omega_\F = 0.
\]
The singular scheme $\Sing(\F)$ of the foliation is the singular scheme of the sheaf $\Omega^1_\F$, which coincides with the vanishing locus of $\omega_\F$.

As mentioned in the introduction, we are interested in extending $\F$ to higher dimensional varieties (not necessarily deformations).

\begin{definition}
Let $\F$ be a codimension $1$ distribution on a smooth closed subvariety $X\subseteq Y$. An \textit{extension $\F^{(n)}$ of order $n$} of $\F$ is given (up to multiplication by a global unit of $\OO_{X^{(n)}}$) by a twisted $1$-form
\[
\omega^{(n)}\in H^0\Big(X^{(n)},\Omega^1_{X^{(n)}}\otimes N_{\F^{(n)}}\Big)
\]
satisfying $i^*\omega^{(n)}=\omega_\F$, where $i:X\hookrightarrow X^{(n)}$ is the natural inclusion. We say that  $\F^{(n)}$ is \textit{integrable} if
$\omega^{(n)}\wedge d\omega^{(n)} = 0.$
\end{definition}  

\begin{remark} \label{def_unf}
If $\F$ is a foliation on $X$, we will mainly be interested in its integrable extensions. In the literature, an extension of $\F$ to a neighborhood of $X$ as a fiber of a deformation $\X\rightarrow T$ is called an \textit{unfolding} of $\F$. 
\end{remark}

\begin{remark} \label{notacion:formasdif}
For computational purposes, it is desirable to have explicit formulas for all possible extensions $\omega^{(n)}$ of a given differential form $\omega^{(n-1)}$ to $X^{(n)}$. Suppose $X\subseteq Y$ is a codimension $1$ regular embedding. Let $\mathcal{U} = \{ U_i \}$ be an affine cover of $X$, and let $t_i \in (I_{X/Y})|_{U_i}$ be local generators as in Remark \ref{los_t's}. Then, on each $U_i$, we have an exact sequence
\[ 0\rightarrow (t_i)\rightarrow \OO_{X^{(n)}}|_{U_i}\rightarrow \OO_X|_{U_i}\rightarrow 0, \]
endowing $\OO_{X^{(n)}}|_{U_i}$ with the structure of a deformation of $U_i$ over $\Spec \mathbb{C}[t_i]/(t_i^{n+1})$. Since each $U_i$ is smooth, these deformations are necessarily trivial, so we can find isomorphisms $X^{(n)}|_{U_i} \simeq U_i\times \Spec \C[t_i]/(t_i^n)$.
However, these isomorphisms are non-canonical, so one would need to carefully keep track of the choices made (cf. Lemma \ref{automorfismos} below). To circumvent this issue, we will instead compare two given extensions $\omega^{(n)}$ and $\overline{\omega}^{(n)}$ of $\omega^{(n-1)}$. Since both coincide when restricted to $X^{(n-1)}$, the conormal exact sequence
\[
0\rightarrow N_{X/Y}^{\otimes - n} \rightarrow \Omega^1_{X^{(n)}}\big|_{X^{(n-1)}}\rightarrow \Omega^1_{X^{(n-1)}}\rightarrow 0
\]
implies that on each $U_i$ we must have 
\begin{equation}
\overline{\omega}^{(n)}_i = \omega^{(n)}_i + \eta_i t_i^n + h_i dt_i^n    
\end{equation}
for some $\eta_i\in \Omega^1_{U_i}$ and $h_i\in \OO_{U_i}$. This notation will be used consistently throughout the paper.
\end{remark}

\begin{lemma} \label{automorfismos}
Let $X^{(n-1)}\subseteq X^{(n)}$ be successive infinitesimal neighborhoods of a codimension $1$ regular embedding $X\subseteq Y$. The space of twisted vector fields \smash{$H^0\big(X,T_X\otimes N_{X/Y}^{\otimes -n}\big)$}
is isomorphic to the group of automorphisms of $X^{(n)}$ whose restriction to $X^{(n-1)}$ is the identity. 
\end{lemma}
\begin{proof}
We will keep the notation of the previous remark. 
Let $\theta$ be an automorphism of $X^{(n)}$ restricting to the identity on $X^{(n-1)}$. From the exact sequence
\[ 0\rightarrow (t_i^n)\rightarrow \OO_{X^{(n)}}\big|_{U_i}\rightarrow \OO_{X^{(n-1)}}\big|_{U_i}\rightarrow 0 \]
and \cite[Lemma 1.2.6]{sernesi2007deformations} we can conclude that 
\begin{equation*}
\theta|_{U_i} = 1+t_i^n v_i
\end{equation*}
for some $v_i\in T_{U_i} = \Der(\OO_{U_i},\OO_{U_i})$. Since these automorphisms coincide on $U_{ij}$ and $t_j^n=f_{ji}^nt_i^n$, we must also have $v_{j}=f_{ij}^nv_i$. This is, $\{v_i\}$ defines a global section of $T_X\otimes N_{X/Y}^{\otimes -n}$.
Reciprocally, using the equation above every such $v=\{v_i \}$ defines an automorphism $\theta_v$ of $X^{(n)}$ relative to $X^{(n-1)}$, establishing a bijection $v\mapsto \theta_v$. From the identity  $\theta_v\theta_{v'}|_{U_i} = (1+t_i^n v_i)(1+t_i^n v_i') = 1 + t_i^n(v_i + v_i') = \theta_{v+v'}|_{U_i}$, we can conclude that 
this map is indeed a group isomorphism.
\end{proof}

\begin{remark} \label{trivialextensions}
Let us briefly observe what the last argument is telling us about our extension problem.
    Using the notation of the above proof, if $\F^{(n)}$ is a distribution on $X^{(n)}$ induced by a twisted differential form $\omega^{(n)}$, then every $\overline{\omega}^{(n)}$ inducing a distribution in the same equivalence class (relative to $X^{(n-1)}$) is locally of the form $\overline{\omega}^{(n)}_i=\omega^{(n)} +t_i^n L_{v_i}\omega + \omega(v_i)dt_i^n$, where $L_{v_i}$ denotes the Lie derivative with respect to the vector field $v_i$. 
\end{remark}

\begin{definition} \label{def:formalext}
    Let $\F$ be a codimension $1$ foliation on a smooth closed subvariety $X\subseteq Y$. A \emph{formal extension} of $\F$ is a a foliation $\smash{\widehat{\F}}$ on $\smash{\widehat{X}}$ extending $\F$, i.e., a sequence of extensions $(X^{(k)},\F^{(k)})$ such that $i_n^*\F^{(n+1)}=\F^{(n)}$, where $i_n:X^{(n)}\hookrightarrow X^{(n+1)}$ is the inclusion of infinitesimal neighborhoods.
\end{definition}

\begin{definition}
 Let $\F$ be a codimension $1$ foliation on a closed subvariety $X\subseteq Y$ and let $k\le n$ be two non-negative integers. Two extensions $\F^{(n)}$ and $\smash{\mybar{1}{0pt}{\F}^{(n)}}$ are \textit{isomorphic relative to $X^{(k)}$} if there is an automorphism $\varphi:X^{(n)}\rightarrow X^{(n)}$ restricting to the identity on $X^{(k)}$ and such that $\smash{\mybar{1}{0pt}{\F}^{(n)}}= \varphi^*\F^{(n)}$. If $k=0$ we will just say that $\F^{(n)}$ and $\smash{\mybar{1}{0pt}{\F}^{(n)}}$ are isomorphic extensions. 
\end{definition}

\subsection{Persistent and Kupka singularities} 

\begin{definition}
Let $\F$ be a codimension $1$ foliation on a smooth variety $X$. A singularity $p\in \Sing(\F)$ is of \textit{Kupka type} if $d\omega_\F(p)\neq 0$. We will denote by $K(\F)$ the set of Kupka points of $\F$.
\end{definition}

Around a Kupka singularity $p\in X$, the foliation can be locally described as the pullback via a smooth morphism $f:(X,p)\rightarrow (\C^2,0)$ of a germ of foliation by curves on the target (see for instance \cite[Theorem 6.13]{neto2020complex}). This is, there exists a coordinate system around $(\C^n,0)\simeq (X,p)$ such that $f$ corresponds to the projection onto the first two variables and $\omega_\F$ takes the form $\omega_\F=A(x_1,x_2)dx_1+B(x_1,x_2)dx_2$. Since $A$ and $B$ vanish simultaneously only at the origin, locally around $p$ we have $K(\F)=f^{-1}(0)$. In particular, the Kupka set is either empty or a (not necessarily connected) locally closed subvariety of pure codimension $2$.

Let $\F$ be a germ of codimension $1$ foliation at $(X,p)$. Recall from Remark \ref{def_unf} that a first order unfolding $\F^{(1)}$ of $\F$ is given (up to multiplication by a unit of $\OO_{(X,p)\times D}$) by an integrable $1$ form $\omega^{(1)}\in\Omega^1_{(X,p)\times D}$ extending $\omega_\F$, where  $D=\Spec(\C[t]/(t^2))$.

\begin{definition} \label{persistent}
A singularity $p$ of a codimension $1$ foliation $\F$ is \emph{persistent} if for every germ of first order unfolding $((X,p)\times D,\F^{(1)}) $, the point $p$ is a singularity of $\F^{(1)}$.
\end{definition}

If $\omega$ is a differential form defining $\F$ at $p$, then following Remark \ref{notacion:formasdif} a first order unfolding will be defined by an integrable differential form $\omega^{(1)}=\omega +t \eta + h dt$ for some $\eta\in\Omega^1_{(X,p)}$ and $h\in\OO_{X,p}$ 

\begin{lemma}\label{eq:unfolding}
    Let $\omega$ be a germ of differentiable $1$-form inducing some foliation $\F$, and let $\eta\in \Omega^1_X$ and $h\in\OO_X$. Then $\omega^{(1)}=\omega+t\eta+hdt\in \Omega^1_{X\times D}$ defines a first order unfolding of $\F$ if and only if 
    $$hd\omega=\omega\wedge(\eta-dh).$$
\end{lemma}
\begin{proof}
First, observe that the form $\omega^{(1)}$ will be integrable if and only if
\begin{align*}
    0&=\omega^{(1)}\wedge d\omega^{(1)}\\
    &=(\omega+t\eta+hdt)\wedge d(\omega+t\eta+hdt)\\
    &=(\eta\wedge d\omega+\omega\wedge d\eta) t + (hd\omega-\omega\wedge(\eta-dh))\wedge dt,
\end{align*}
which holds exactly when both $\alpha:=\eta\wedge d\omega+\omega\wedge d\eta$ and $\beta:=hd\omega-\omega\wedge(\eta-dh)$ vanish identically. For the other implication, let us first analyze the case where $h=0$. In this situation, the equation $\beta=0$ reduces to $\omega\wedge\eta=0$ and hence $\eta=\lambda \omega$ for some $\lambda\in\OO_X^*$. But this means that $\omega^{(1)}=(1+\lambda t)\omega$, which is of course integrable. If on the other hand $h\neq 0$, it is enough to check that these elements satisfy the equation
$$ d\beta+\frac{2}{h}\beta\wedge(\eta-dh)=\alpha,$$
and therefore the vanishing of $\beta$ implies the integrability of $\omega^{(1)}$.
\end{proof}

In terms of the differential forms defining the foliation, the persistency of a foliation singularity $p$ should be interpreted as follows. 
Since $\omega(p)=0$, the point $p$ is said to be a persistent singularity if for every first order unfolding $\omega^{(1)}=\omega+t \eta+hdt$ of $\omega$ we have $\omega^{(1)}(p)=h(p)dt=0$.

\begin{definition}\label{def:Iper}
    The scheme $\Per(\F)$ of persistent singularities of a codimension $1$ foliation $\F$ on $X$ is the subscheme defined by the sheaf of ideals 
    $\I_{\Per}(U)=\{h\in \OO_X(U) : hd\omega_\F=0\in \left(\Omega^2_X/( \omega_\F \wedge \Omega^1_X) \right)(U)\}.$
\end{definition}

The above definition makes sense due to the following statement.

\begin{proposition}\cite[Proposition 3.6]{persistent2021}
    Let $p$ be a point in $\Sing(\F)$. Then $p$ is a persistent singularity if and only if $p \in \Per(\F)$.
\end{proposition}

\begin{remark}
    If $\Sing(\F)$ has codimension greater than two around $p$, then by Malgrange's Theorem (locally around $p$) $\F$ admits a holomorphic first integral. This is, $\F$ is defined by an element of the form $\omega=gdf$ for some $g\in\OO_{X,p}^*$ and some $f\in \OO_{X,p}$. In particular, $p$ is not a persistent singularity, since in this case one can smooth out the singularity by considering a local unfolding of the form $\omega^{(1)}=g d(f+t)\in\Omega^1_{(X,p)\times D}$. We will analyze this situation in more detail in Example \ref{unf_int1}.
\end{remark}

Typically, one should expect a codimension $1$ foliation $\F$ to have persistent singularities. This is the case if $X$ is a smooth projective variety such that \smash{$H^1(X,N_\F^\vee)=0$} and $c_1(N_\F) \neq 0$ \cite[Proposition 3.8]{persistent2021}. The same authors proved that persistent and Kupka singularities are intimately related:

\begin{theorem}\cite[Theorem 3.14]{persistent2021} \label{teo:persitent=kupka}Let $X$ be a smooth projective variety and $\F$ a codimension $1$ foliation on $X$ with reduced singular scheme. Then $\smash{\Per(\F)_{\operatorname{red}}=\overline{K(\F)}}$.
\end{theorem}

\subsection{Sheaf of principal parts along a distribution}
We will end this section by introducing a central tool in our study.
Let $X$ be a smooth variety of dimension $n$ over $\C$. The Chern map 
$$c:H^1(X,\OO_X^*)\to H^1(X,\Omega^1_X)=\Ext^1(T_X,\OO_X)$$ 
is induced at the level of cocycles by the assignment
$$ \{f_{ij}\} \mapsto \bigg\{\frac {df_{ij}}{f_{ij}} \bigg\}.  $$
This is, every line bundle $\L$ on $X$ yields an extension 
$$0\rightarrow \OO_X \rightarrow \E_\L \rightarrow T_X \rightarrow 0, $$
called the \emph{Atiyah extension of} $\L$. The sheaf $\E_\L$ is locally free of rank $n+1$ and 
\[
\PP_\L:=\E_\L^\vee\otimes \L
\]
is called the \emph{sheaf of first-order principal parts of $\L$}. 

Given a codimension $1$ distribution $\F$ and a line bundle $\L$ on $X$, the \textit{sheaf of principal parts of $\L$ along $\F$} is the extension 
\[
0 \rightarrow  \Omega^1_{\F}\big(N_\F\otimes \L\big) \rightarrow {\PP_{\F,\L}} \rightarrow N_\F\otimes \L \rightarrow 0
\]
attached to the image of the first Chern class $c_1(\L)$ under the map
\[
\Ext^1\!\!\big(\OO_X,\Omega^1_X\big)\rightarrow \Ext^1\!\!\big(\OO_X, \Omega^1_\F\big) = \Ext^1\!\!\big(N_\F\otimes \L, \Omega^1_\F(N_\F\otimes \L)\big).\\
\]

\begin{lemma}\label{cohomology_principal_parts}
Given a cover $\U=\{U_i\}$ of $X$ by affine open subsets, the cohomology groups of the sheaf of principal parts of $\L$ along $\F$ can be computed as
{\smaller\[
H^{k}\big(X,\PP_{\F,\L}\big) = \frac{\big\{(\alpha,\beta)\in Z^k(\U,N_\F\otimes \L)\times C^{k}(\U,\Omega^1_\F(N_\F\otimes \L)):\; \delta(\beta) = \big[c_1(\L)\cup\alpha\big]  \big\}}{\big\{\big(\delta(\gamma),[c_1(\L)\cup\gamma]+\delta(\sigma)\big):\,(\gamma,\sigma)\in C^{k-1}(\U,(N_\F\otimes \L)\oplus\Omega^1_\F(N_\F\otimes \L))\big\}},
\]}
where $(C^\bullet(\U,-),\delta)$ are the $\check{C}$ech complexes associated to $\U$.
\end{lemma}

\begin{proof}
Take transition functions $\{f_{ij}\}$ of the line bundle $\L$ with respect to $\U$. By an abuse of notation, we will identify $c_1(\L)\in H^1(X,\Omega^1_X)$ with its $\check{C}$ech 1-cocycle representative $\{df_{ij}/f_{ij}\}\in Z^{1}(\U,\Omega^1_X)\subseteq C^{1}(\U,\Omega^1_X)$. The image of the class $c_1(\L)$ in the group $\Ext^1\!\!\big(N_\F\otimes \L, \Omega^1_\F(N_\F\otimes \L)\big)$ is identified with the morphism (in the derived category of abelian groups $D(\operatorname{Ab})$) given by the cup product
\[
\big[c_1(\L)\cup -\big]: C^{\bullet}\big(\U,N_\F\otimes \L\big) \rightarrow  C^\bullet\big(\U,\Omega^1_\F(N_\F\otimes \L)\big) \big[1\big].
\]
See \cite[Prop. 1.4.3]{dimca2004sheaves} for details. By definition of the sheaf of principal parts,
\begin{equation*}\label{derived}
C^{\bullet}\big(\U,N_\F\otimes \L\big) \xrightarrow[]{c_1(\L)} C^{\bullet}\big(\U,\Omega^1_\F(N_\F\otimes \L)\big)\big[1\big]\rightarrow C^{\bullet}\big(\U,\PP_{\F,\L}\big)\big[1\big]\rightarrow  C^{\bullet}\big(\U,N_\F\otimes \L\big) \big[2\big]
\end{equation*}
is a triangle in $D(\operatorname{Ab})$. Recall that the \textit{cone} of a morphism $f: A^\bullet \rightarrow B^\bullet$ between chain complexes in a given abelian category is the chain complex
\[
\operatorname{Cone}^\bullet(f) = A^\bullet[1] \oplus B^\bullet
\]
with differential $\delta_f(a,b) = (-\delta(a),\delta(b)-f(a))$\footnote{The literature may adopt a different sign convention for the differential in the mapping cone, defining $\operatorname{Cone}^\bullet(f)$ instead as $C^\bullet(-f)$. However, this distinction is non-essential: the two complexes are isomorphic via the explicit map $(a,b)\mapsto (-a,b)$ \cite[p. 7]{dimca2004sheaves}}. If we take $f = \big[c_1(\L)\cup -\big]$ as before, the above triangle implies that the shifted cone $\operatorname{Cone}^\bullet(f)[-1]$ is isomorphic to $C^\bullet(\U,\PP_{\F,\L})$ in $D(\operatorname{Ab})$. In particular, these complexes have the same cohomology, which can be computed as stated in the Lemma.
\end{proof}

\section{Infinitesimal extensions of distributions and foliations}\label{section3}
We will now develop a key ingredient of our strategy, namely the understanding of the extension problem at the level of the infinitesimal neighborhoods of $X$. Throughout this section we will always assume that $X\subset Y$ is a regularly embedded divisor. Just as in the previous section, no compactness assumptions are needed at this stage. 

\subsection{Extensions and obstructions}
In this section we will compute the first order extensions of a foliation $\F$ and develop an obstruction theory for extending to other infinitesimal neighborhoods of $X$. 

\begin{proposition} \label{extensions:dist}
Let $X^{(n-1)}\subseteq X^{(n)}$ be successive infinitesimal neighborhoods of $X\subseteq Y$. If $\F^{(n-1)}$ is an extension of a codimension $1$ distribution $\F$ on $X$, then there is a free transitive action of the vector space 
\[
H^0\big(X,\PP_{\F,N_{X/Y}^{\otimes - n}}\big)
\]
on the set of extensions of $\F^{(n-1)}$ along $X^{(n)}$. 
\end{proposition}

\begin{proof}
Let $\omega$ be a global section of $\Omega^1_X\otimes N_\F$ defining $\F$ and $\U = \{ U_i\}$ be an affine cover of $X$. Let $\{g_{ij}\}$ and $\{f_{ij} \}$ be cocycles representing $\L^{(0)}:=N_\F$ and $N_{X/Y}^{\otimes -1}$ respectively and $\omega^{(n)},\overline{\omega}^{(n)}$ be two twisted differential one forms on $X^{(n)}$ extending $\omega$. Since these differential forms coincide after pulling back to $X^{(n-1)}$, then by Remark \ref{notacion:formasdif} on each $U_i$ we must have 
\begin{equation}\label{eq:transitividad}
\overline{\omega}^{(n)}_i = \omega^{(n)}_i + \eta_i t_i^n + h_i dt_i^n    
\end{equation}
for some $\eta_i\in \Omega^1_{U_i}$ and $h_i\in \OO_{U_i}$. We will now see that $\{h_i,[\eta_i]\}$ defines a global section of the corresponding sheaf of principal parts. If $\omega^{(n)},\overline{\omega}^{(n)}$ are twisted by line bundles \smash{$\L^{(n)}$, $\mybar{1}{0pt}{\L}^{(n)}$} with transition functions $\big\{g_{ij}^{(n)}\big\}$, $\big\{\overline{g}_{ij}^{(n)}\big\}$ respectively, there should $1$-cocycle $\{u_{ij}\}$ satisfying the equations (\ref{eq:1}). Expanding the identities
\begin{equation*}
\overline{\omega}^{(n)}_j-\omega_j^{(n)}= \overline{g}_{ij}^{(n)}\overline{\omega}^{(n)}_i-g_{ij}^{(n)}\omega^{(n)}_i
\end{equation*}
we arrive at the relations
\begin{align*}
\eta_j t_j^n + h_j dt_j^n & = g_{ij}^{(n)}(1+u_{ij} t_i^n)\big(\omega^{(n)}_i + \eta_i t_i^n + h_i dt_i^n\big)-g_{ij}^{(n)}\omega^{(n)}_i\\
 & = g_{ij}^{(n)}\Big[\big(\eta_i t_i^n + h_i dt_i^n\big)+ u_{ij} t_i^n\big(\omega^{(n)}_i + \eta_i t_i^n + h_i dt_i^n\big)\Big]\\
 & = g_{ij}\;\Big[\big(\eta_i t_i^n + h_i dt_i^n\big)+ u_{ij}\omega_it_i^n\Big]\\
 & = g_{ij}\;\Big[\big(\eta_i+u_{ij}\omega_i\big) t_i^n + h_i dt_i^n\Big],
\end{align*}
where in the second to last equality we used that $I_{X/Y}^{n+1}=0$ on $X^{(n)}$ and both $g_{ij}^{(n)} t_i^n = g_{ij} t_i^n$ and $g_{ij}^{(n)}dt^n_i=g_{ij} dt_i^n$.
Now since $t_i^n= f_{ij}^n t_j^n$, 
\begin{align*}
\eta_j t_j^n + h_j dt_j^n& = g_{ij}\Big[\big(\eta_i+u_{ij}\omega_i\big)f_{ij}^n t_j^n + h_i\big(f_{ij}^ndt_j^n+t_j^ndf_{ij}^n\big)\Big]\\
     &= \Big(
     g_{ij} f_{ij}^n\big(\eta_i+u_{ij}\omega_i\big)+h_i g_{ij} df_{ij}^n\Big)t_j^n + \Big(g_{ij} f_{ij}^nh_i\Big)dt_j^n.
\end{align*}
From this we can conclude that the local sections $(h_i,\eta_i)$ satisfy the cocycle conditions $h_j = f_{ij}^ng_{ij}h_i$ and 
\[ 
\eta_j - f_{ij}^ng_{ij}\big(\eta_i+u_{ij}\omega_i\big) = h_i g_{ij} f_{ij}^n\frac{d f_{ij}^n}{f_{ij}^n} = \bigg\{\frac{d f_{ij}^n}{f_{ij}^n}\bigg\}\cup \{h_j\}.
\]
If we denote by $[\sigma]$ the class of a local 1-form $\sigma$ in $\Omega^1_{\F}$, the previous sections satisfy 
\[
\delta\big(\{h_i\}\big) = 0,\quad\text{and}\quad\delta\big(\{[\eta_i]\}\big) = \big[c_1(N_{X/Y}^{\otimes - n})\cup\{h_i\}\big], 
\]
where $\delta$ denotes, making a notional abuse, both the $\check{C}$ech differentials of the sheaves $\L\otimes N_{X/Y}^{-1}$ and $\Omega^{1}_\F(\L\otimes N_{X/Y}^{\otimes - n})$. In conclusion, the collection ${(h_i,[\eta_i])}$ defines a global section of $\PP_{\F,N_{X/Y}^{\otimes - n}}$ by the Lemma \ref{cohomology_principal_parts} in the previous section. Reading the above argument backwards we see that for any such cocycle $\{(h_i,[\eta_i])\}$ and an extension $\omega^{(n)}$, the cocycle $\big\{\omega_i^{(n)}+\eta_it_i^n + h_i dt_i^n\big\}$ defines a differential $1$-form twisted by the line bundle defined with transitions functions $\big\{g_{ij}^{(n)}(1+u_{ij}t_i^n)\big\}$ which also extends $\omega^{(n-1)}$. On the other hand, acting on an extension $\omega^{(n)}$ by a cocycle of the form $\{(0,s_i\omega_i)\}$ does not change the corresponding distribution on $X^{(n)}$ because 
\[
\omega_i^{(n)}+s_i\omega_it_i^n = (1+s_it_i^n)\omega_i^{(n)},
\]
and $1+s_it_i^n$ is a unit on $X^{(n)}$. This implies that \smash{$H^0\big(X,\PP_{\F,N_{X/Y}^{\otimes - n}}\big)$} acts on the set of extensions of $\F^{(n-1)}$. This action is transitive by equation (\ref{eq:transitividad}) and the argument above. Finally, we claim this action is also free. Indeed, if some global section $\{(h_i,[\eta_i])\}$ fixes some extension $\F^{(n)}$ (induced by $\omega^{(n)}$), we should have 
$$\lambda_i\omega_i^{(n)} = \omega_i^{(n)}+\eta_it_i^n + h_i dt_i^n$$
on each $U_i$ where $\lambda_i$ is a local unit on $X^{(n)}$. Since both sides of the equation extend $\omega^{(n-1)}$, the unit $\lambda_i$ must be of the form $1+s_it_i^{n}$ for some $s_i\in \OO_{U_i}$ 
 and therefore $(h_i,[\eta_i]) = (0,[s_i\omega_i])= 0$ in $\PP_{\F,N_{X/Y}^{\otimes - n}}$.
\end{proof}

Now that we understand the infinitesimal extensions of successive order of a given distribution, the next step would be to identify the ones that are in the same equivalence class.

Recall from Lemma \ref{automorfismos} that the automorphism group of $X^{(n)}$ relative to $X^{(n-1)}$ is parametrized by $H^0(X,T_X\otimes N_{X/Y}^{\otimes -n})$. In particular, given an extension $\omega^{(n)}$ of $\omega^{(n-1)}$ and a section $v\in H^0(X,T_X\otimes N_{X/Y}^{\otimes -n})$ one can consider the extension $\theta_v^* \omega^{(n)}$, where $\theta_v:X^{(n)}\to X^{(n)}$ is the automorphism induced by $v$. By Remark \ref{trivialextensions}, the difference between these extensions will coincide with the image of $v$ under the $\C$-linear map
\begin{align*}
    H^0(X,T_X\otimes N_{X/Y}^{\otimes -n})&\to H^0(X,\PP_{\F,N_{X/Y}^{\otimes - n}})\\
    v&\mapsto \big(\omega_i(v_i),[L_{v_i}\omega_i]\big),
\end{align*}
where $v_i$ and $\omega_i$ are the restrictions of $v$ and $\omega$ to an affine cover $\U = \{ U_i\}$.
This discussion shows the following.
\begin{proposition}
    There is a 1-1 correspondence between equivalence classes of extensions of $\F^{(n-1)}$ to $X^{(n)}$ and elements in the quotient
$H^0\big(X,\PP_{\F,N_{X/Y}^{\otimes - n}}\big)/H^0\big(X,T_X\otimes N_{X/Y}^{\otimes -n}\big)$.
\end{proposition}

In order to determine which of these extensions correspond to integrable extensions, let us first compute the elements in $H^0(X,\PP_{\F,N_{X/Y}^{\otimes - n}})$ that preserve integrability.

\begin{lemma} \label{eq:integrabilidad}
    Let $\omega$ be a local integrable $1$-form and $\omega^{(n)}$ a local integrable extension. Then $\overline{\omega}^{(n)}=\omega^{(n)} + t^{n} \eta + h dt^n$ is integrable if and only if 
    $$ h d\omega = \omega \wedge ( \eta-dh).$$
    Moreover, the choice of $h$ completely determines the class $[\eta]\in \Omega^1_\F$.
\end{lemma}
\begin{proof}
Observe that 
\begin{align*}\overline{\omega}^{(n)}\wedge d\overline{\omega}^{(n)} &= \big(\omega^{(n)} + t^n \eta + h dt^n\big)\wedge \big(d\omega^{(n)} + dt^n\wedge \eta + t^n d\eta + dh\wedge dt^n\big)  \\
&= \Big( \eta\wedge d\omega^{(n)} + \omega^{(n)}\wedge d\eta\Big)t^n + \Big(h d\omega^{(n)} - \omega^{(n)}\wedge (\eta-dh)\Big)\wedge dt^n. 
\end{align*}
Since $t^{n+1}=0$, this actually reduces to the desired expression, from where we can proceed in the same manner as in Lemma \ref{eq:unfolding}. The uniqueness of $[\eta]$ can be checked as in \cite{suwa1992unfoldings}.
\end{proof}
The lemma suggests that the subspace of $H^0(X,\PP_{\F,N_{X/Y}^{\otimes - n}})$ preserving integrability on extensions of $\F$ corresponds to the kernel of the $\C$-linear map 
\begin{align*}
   \PP_{\F,N_{X/Y}^{\otimes -n}} & \longrightarrow \Omega^2_X(N_\F^{\otimes 2}\otimes N_{X/Y}^{\otimes -n}) \\
    ([\eta],h) &\longmapsto h d\omega - \omega \wedge ( \eta-dh).
\end{align*}
In fact, one can further compose with the projection 
\begin{equation}\label{def:QF}
\Omega^2_X(N_\F^2\otimes N_{X/Y}^{\otimes - n}) \to \Omega^2_X(N_\F^{\otimes 2}\otimes N_{X/Y}^{\otimes - n}) / \big(\omega_\F\wedge \Omega^1_X(N_\F\otimes N_{X/Y}^{\otimes - n})\big)=:Q_{\F}(n) 
\end{equation}
in order to end up with a $\OO_X$-linear morphism
\[
    d\omega_{N_{X/Y}^{\otimes - n}}: \PP_{\F,N_{X/Y}^{\otimes - n}}  \longrightarrow Q_\F(n)
\]
defined at the level of local sections by the rule $([\eta],h)\mapsto [h d\omega]$.
We will use the notation $\PP_{\F,N_{X/Y}^{\otimes - n}}^{\,\Int}$ for the kernel of this map. 

\begin{proposition} \label{def:fol}Let $X^{(n-1)}\subseteq X^{(n)}$ be be successive infinitesimal neighborhoods of $X\subseteq Y$. Suppose $\F^{(n-1)}$ is an extension on $X^{(n-1)}$ of a codimension $1$ foliation $\F$. Then there is a free and transitive action of
\[
H^0\big(X,\PP_{\F,N_{X/Y}^{\otimes - n}}^{\,\Int}\big)
\]
on the set of integrable extensions of $\F^{(n-1)}$ to $X^{(n)}$.
Moreover, there is a $1$-$1$ correspondence between the set of equivalence classes of integrable extensions and the vector space 
$$\Ex^n(\F,N_{X/Y}) := \frac{H^0\Big(X,\PP_{\F,N_{X/Y}^{\otimes - n}}^{\,\Int}\Big)}{H^0\big(X,T_X\otimes N_{X/Y}^{\otimes -n}\big)}.$$ 
\end{proposition}
\begin{proof}
This is a direct application of Proposition \ref{extensions:dist} and Lemma \ref{eq:integrabilidad}.
\end{proof}

\begin{remark}
    In the case where $X$ is a fiber of a submersion $Y\rightarrow C$ onto a smooth curve, we have $N_{X/Y}\simeq \OO_X$ and therefore all the groups $\Ex^n(\F,N_{X/Y})$ coincide. We will hence denote them by $\Ex(\F)$. 
\end{remark}

Let us now move on to computing the obstructions to extend along a small extension of infinitesimal thickenings.

\begin{proposition} \label{obst:dist}
Let $X^{(n-1)}\subseteq X^{(n)}$ be successive infinitesimal neighborhoods of $X\subseteq Y$. The obstruction to extend a distribution $\F^{(n-1)}$ on $X^{(n-1)}$ to $X^{(n)}$ lies in the group
\[
H^1\Big(X, \PP_{\F,N_{X/Y}^{\otimes - n}}\Big).
\]
\end{proposition}
\begin{proof} Suppose $\omega^{(n-1)}$ is a $1$-form twisted by $\L^{(n-1)}$ defining the distribution $\F^{(n-1)}$, and $\big\{g_{ij}^{(n-1)}\big\}$ are the transition functions of this line bundle. Choose extensions $\omega_i^{(n)}$ of the local $1$-forms $\omega_i^{(n-1)}$ to $X^{(n-1)}$, and units $g_{ij}^{(n)}$ extending $g_{ij}^{(n-1)}$. If we were lucky enough, this data would glue together, promoting both a line bundle extending $\L^{(n-1)}$ and a 1-form twisted by this line bundle and extending $\omega^{(n-1)}$. In any case, since these local sections glue when restricted to $X^{(n-1)}$, by Remark \ref{notacion:formasdif} over each $U_{ij}$ they must satisfy relations of the form
\begin{equation}\label{Obstructions:eq1}
\omega^{(n)}_j - g_{ij}^{(n)}\omega_i^{(n)} = \eta_{ij} t_j^n+ h_{ij}dt_j^n,
\end{equation}
for some local sections $\eta_{ij}\in \Omega^1_{U_{ij}}$ and $h_{ij}\in \OO_{U_{ij}}$.
Being extensions of $\{g_{ij}^{(n-1)} \}$, on $U_{ijk}$ the new transition functions must also satisfy
\begin{equation}\label{Obstructions:eq2}
g_{ij}^{(n)}g_{jk}^{(n)}g_{ki}^{(n)} = 1 + u_{ijk}t_k^n,
\end{equation}
for some $u_{ijk}\in\OO_{U_{ijk}}$.
We claim the collection \smash{$\{(h_{ij},[\eta_{ij}])\}$} is a \smash{$\check{C}$}ech $1$-cocycle of the sheaf \smash{$\PP_{\F,N_{X/Y}^{\otimes - n}}$}. In order to prove this, let us rearrange the expression 
\[
g_{jk}^{(n)}\big(\omega_j^{(n)}-g_{ij}^{(n)}\omega_i^{(n)}\big) - \big(\omega_k^{(n)}-g_{ik}^{(n)}\omega_i^{(n)}\big) + \big(\omega_k^{(n)}-g_{ij}^{(n)}\omega_j^{(n)}\big)
\]
in two different ways. On one hand, some of the terms cancel each other out and the only ones remaining are
\[
\Big(g_{ik}^{(n)} - g_{ij}^{(n)}g_{jk}^{(n)}\Big) \omega_i^{(n)} = \Big(-u_{ijk}g_{ik}\omega_i\Big)t_k^n.
\]
On the other hand, appealing to the defining equations of the local sections $h_{ij}$ and $\eta_{ij}$, we obtain the expression
\[
g_{jk}\big(\eta_{ij} t_j^n+ h_{ij}dt_j^n\big) - \big(\eta_{ik} t_k^n + h_{ik}dt_k^n\big) + \big(\eta_{jk} t_k^n + h_{jk}dt_k^n\big).
\]
Rearranging the terms and using the fact that $t_j^n=f_{jk}^nt_k^n$, this becomes
\[
\Bigg(f_{jk}^ng_{jk}\eta_{ij}+f_{jk}^ng_{jk}h_{ij}\frac{df_{jk}^n}{f_{jk}^n}-\eta_{ik}+\eta_{jk}\Bigg) t_k^n + \Big(g_{jk}h_{ij}-h_{ik}+h_{jk}\Big) dt_k^n.
\]
Comparing the coefficients of $t_k^n$ and $dt_k^n$ in each result we can conclude that
\begin{equation*}
\delta\big(\{h_{ij}\}\big) = 0 \quad\text{and}\quad \delta\big(\{[\eta_{ij}]\}\big) = \big[c_1(N_{X/Y}^{\otimes - n})\cup \{h_{ij}\}\big]
\end{equation*} 
as promised. We claim that the class of this $1$-cocycle in \smash{$H^1\big(X,\PP_{\F,N_{X/Y}^{\otimes - n}}\big)$} does not depend of the choices made, only on the distribution $\F^{(n-1)}$ itself. As seen along the proof of Proposition \ref{extensions:dist}, alternative extensions \smash{$\overline{\omega}_i^{(n)}$} and \smash{$\overline{g}_{ij}^{(n)}$} to \smash{$\omega_i^{(n-1)}$} and \smash{$g_{ij}^{(n-1)}$} respectively must be written as
\begin{equation}\label{Obstructions:eq3}
\overline{\omega}_i^{(n)} = \omega_i^{(n)}+ \eta_it_i^n+h_idt_i^n\quad\text{and}\quad \overline{g}_{ij}^{(n)} = g_{ij}^{(n)}\big(1+u_{ij}t_i^n\big) 
\end{equation}
for some $\eta_i\in \Omega^1_{U_i}$ and $u_{ij}\in \OO_{U_{ij}}$.
In the same fashion as in (\ref{Obstructions:eq1}), these choices also verify compatibility equations
\begin{equation}\label{Obstructions:eq4}
\overline{\omega}^{(n)}_j - \overline{g}_{ij}^{(n)}\overline{\omega}_i^{(n)} = \overline{\eta}_{ij} t_j^n+ \overline{h}_{ij}dt_j^n.  
\end{equation}
for some local sections $\overline{\eta}_{ij}$ and $\overline{h}_{ij}$. Combining equations (\ref{Obstructions:eq1}), (\ref{Obstructions:eq3}) and (\ref{Obstructions:eq4}), and comparing the coefficients of $t_j^n$ and $dt_j^n$ we get the following relations 
\[
\overline{h}_{ij}-h_{ij} = h_{j}-f_{ij}^ng_{ij}h_i\quad\text{and}\quad \overline{\eta}_{ij}-\eta_{ij} = \eta_{j}-f_{ij}^ng_{ij}\big(\eta_{i}+u_{ij}\omega_i\big) -h_i f_{ij}^ng_{ij} \frac{df_{ij}^n}{f_{ij}^n},
\]
which in turn imply
\[
\overline{h}_{ij}-h_{ij} =\delta\big(\{h_{i}\}\big)\quad\text{and}\quad \big[\overline{\eta}_{ij}\big]-\big[\eta_{ij}\big]=\delta\big(\{[\eta_{i}]\}\big)+\big[c_1(N_{X/Y}^{\otimes - n})\cup \{h_{i}\}\big].
\]
Therefore, following Lemma \ref{cohomology_principal_parts} once again, the $1$-cocycles \smash{$\{(h_{ij},[\eta_{ij}])\}$} and  \smash{$\{(\overline{h}_{ij},[\overline{\eta}_{ij}])\}$} represent the same element in \smash{$H^1\big(X,\PP_{\F,N_{X/Y}^{\otimes - n}}\big)$} as expected. In particular, this implies that if the distribution $\F^{(n-1)}$ can be extended, then this obstruction must be zero. Reciprocally, if the obstruction \smash{$\{(h_{ij},[\eta_{ij}])\}$} is zero, we can find local sections $h_i$, $\eta_i$, $u_{ij}$ such that
\[
-h_{ij} = h_{j}-f_{ij}^ng_{ij}h_i\quad\text{and}\quad -\eta_{ij} = \eta_{j}-f_{ij}^ng_{ij}\big(\eta_{i}+u_{ij}\omega_i\big) -h_i f_{ij}^ng_{ij} \frac{df_{ij}^n}{f_{ij}^n},
\]
and consequently the modified local sections
\[
\omega_i^{(n)}+ \eta_it_i^n+h_idt_i^n\quad\text{and}\quad g_{ij}^{(n)}\big(1+u_{ij}t_i^n\big) 
\]
glue together.
\end{proof}

\begin{remark} One of the most elementary reasons why a distribution may not extend is the non-extendability of the normal bundle $\L = N_{\F}$. With the same approach of the proof of Proposition \ref{ext_line_bundles}, one can show that the vector space $H^2(X,N_{X/Y}^{\otimes - n})$ is an obstruction space for this problem (this is also stated in \cite[Section 5.a]{griffiths1966extension}). These obstructions are related in the following manner. Consider the exact sequence
\[
0\rightarrow N_{X/Y}^{\otimes - n}\xrightarrow{\omega_\F} \PP_{N_{X/Y}^{\otimes - n}}\otimes N_\F\rightarrow \PP_{\F,N_{X/Y}^{\otimes - n}}\rightarrow 0.
\]
Looking at the long exact sequence in cohomology we get a map 
\[
H^0\big(X,\PP_{\F,N_{X/Y}^{\otimes - n}}\big)\rightarrow H^1(X,N_{X/Y}^{\otimes - n})
\]
sending the extension $\F^{(n)}$ to its associated extension of the line bundle $\L^{(n-1)}$ to $X^{(n)}$ and a map 
\[
H^1\big(X,\PP_{\F,N_{X/Y}^{\otimes - n}}\big)\rightarrow H^2(X,N_{X/Y}^{\otimes - n})
\]
pointing at the obstruction to extend $\L^{(n-1)}$. 
\end{remark}

\begin{example}\cite[Section 2.1]{cerveau2013quelques}\label{ej:cerveau} 
Consider the foliation $\F$ on $X = \P^1\times \P^1$ given by one of the rulings. Each line bundle on $\P^1\times\P^1$ is of the form $\OO(k_1,k_2)=\OO(k_1)\boxtimes\OO(k_2)$, and the normal bundle of $\F$ can be identified with $\OO_X(2,0)$. This surface can be embedded in $Y = \P^3$ as a smooth quadric with normal bundle $N_{X/Y} = \OO(1,1)$. Since $H^2(X,\OO_X(-n,-n))=0$ for $n>1$, we see that the first (and only) possible obstruction for extending line bundles appears in the first infinitesimal neighborhood.
On the other hand, the pullback of line bundles maps $\OO_Y(1)$ to $\OO_X(1,1)$, and therefore $N_\F$ is not the restriction of a line bundle on $Y$. This implies that $\F$ does not extend to $Y$ - not even extend to $X^{(1)}$, since by the remark above the corresponding obstruction in $H^1(X,\PP_{\F,N_{X/Y}^{\otimes - 1}})$ maps to a non-zero element in $H^2(X,\OO_X(-1,-1))$.
\end{example}

Observe that a key point in the proof of the Proposition above is the fact that the extension problem for codimension $1$ \emph{distributions} is locally unobstructed. This is not the case in the context of foliations. 

\begin{definition} \label{def:unobstruced}
    Let $\F$ be a germ of foliation on $X=(\C^m,0)\subseteq (\C^{m+1},0)$. We say that $\F$ is \textit{unobstructed} if for every $n$ and every extension $\F^{(n-1)}$ of order $n-1$ of $\F$, there exists an extension $\F^{(n)}$ on $X^{(n)}$ restricting to $\F^{(n-1)}$ on $X^{(n-1)}$. 
\end{definition}

\begin{remark}
    Although they are closely related, the above definition does not correspond to the unobstructedness of the deformation functor of the germ of foliation $\F$. For a first approach to the comparison of these problems the reader is referred to \cite[Section I.7]{suwa1992unfoldings}.
\end{remark}

\begin{proposition}\label{obst:h1_pint}
Suppose $\F$ is a foliation on $X$ with unobstructed singularities. The obstruction to further extend an extension $\F^{(n-1)}$ on $X^{(n-1)}$ to  $X^{(n)}$ lies in the group 
    $$H^1\big(X,\PP_{\F,N_{X/Y}^{\otimes - n}}^{\,\Int}\big).$$ 
\end{proposition}
\begin{proof}
We will follow the strategy used for Proposition \ref{obst:dist}.  Suppose $\F^{(n-1)}$ is an extension induced by $\omega^{(n-1)}\in H^0(X^{(n-1)},\Omega^1_{X^{(n-1)}}\otimes \L^{(n-1)})$ and $\big\{g_{ij}^{(n-1)}\big\}$ are the transition functions for $\L^{(n-1)}$ with respect to same affine cover $\U$. Let $\omega_i^{(n)}$ be local \emph{integrable} extensions of $\omega_i^{(n-1)}$ and units $g_{ij}^{(n)}$ extending $g_{ij}^{(n-1)}$. According to the calculations above, these local extensions will induce an element $[(h_{ij},[\eta_{ij}])]\in H^1\big(X,\PP_{\F,N_{X/Y}^{\otimes - n}}\big)$,
where 
$$\omega^{(n)}_j = g_{ij}^{(n)}\omega_i^{(n)} + \eta_{ij} t_j^n+ h_{ij}dt_j^n. $$
Observe that since both $\omega^{(n)}_j$ and $g_{ij}^{(n)}\omega_i^{(n)}$ are integrable, by Lemma \ref{eq:integrabilidad} we also know that $h_{ij}d\omega=0\in \Omega_X^2/(\omega_\F\wedge \Omega^1_X)$, implying that these cocycles also define a class in $H^1(X,\PP_{\F,N_{X/Y}^{\otimes - n}}^{\,\Int})$. From here we can conclude just as in the previous argument, keeping in mind that the local corrections $\omega_i^{(n-1)}+ \eta_it_i^n+h_idt_i^n$ will turn out to be integrable for any family of local sections $(\eta_i,h_i)\in \PP_{\F,N_{X/Y}^{\otimes - n}}^{\,\Int}$.
\end{proof}

\subsection{The role of persistent singularities}

 A careful reader may have noticed that the definition of the morphism $d\omega_{N_{X/Y}^{\otimes - n}}: \PP_{\F,N_{X/Y}^{\otimes - n}} \to  Q_{\F}(n) $ is independent of the $\eta$-coordinate of the elements in $\PP_{\F,N_{X/Y}^{\otimes - n}}$. Indeed, this morphism factorizes through the projection to the h-coordinate \smash{$\PP_{\F,N_{X/Y}^{\otimes - n}}\to N_\F\otimes N_{X/Y}^{\otimes - n}$}.

 \begin{proposition}\label{pint=per} The projection $\PP_{\F,N_{X/Y}^{\otimes - n}} \to N_\F\otimes N_{X/Y}^{\otimes - n}$ induces an isomorphism 
 $$\PP_{\F,N_{X/Y}^{\otimes - n}}^{\,\Int}\simeq \I_{\Per} \otimes N_\F\otimes N_{X/Y}^{\otimes - n},$$
    where $\I_{\Per}$ is the ideal of persistent singularities.  
 \end{proposition}
\begin{proof}
The aforementioned morphism is the dashed arrow in the commutative diagram with exact rows
\begin{center}
\begin{tikzcd}
0 \arrow[r] & {\PP_{\F,N_{X/Y}^{\otimes - n}}^{\,\Int}} \arrow[r] \arrow[d, dashed] & {\PP_{\F,N_{X/Y}^{\otimes - n}}} \arrow[r, "d\omega_\F"] \arrow[d] & Q_{\F}(n) \arrow[d,equal] \\
0 \arrow[r] & \I_{\Per} \otimes N_\F\otimes N_{X/Y}^{\otimes - n} \arrow[r]         & N_\F\otimes N_{X/Y}^{\otimes - n} \arrow[r, "d\omega_\F"]          & Q_{\F}(n),
\end{tikzcd}
\end{center}
whose exactness of follow from the definition on $\PP_{\F,N_{X/Y}^{\otimes - n}}^{\,\Int}$ and Definition \ref{def:Iper} respectively.
The uniqueness of $[\eta]$ in Lemma \ref{eq:integrabilidad} automatically implies that the projection $$\PP_{\F,N_{X/Y}^{\otimes - n}}^{\,\Int}\to \I_{\Per} \otimes N_\F\otimes N_{X/Y}^{\otimes - n}$$
is a monomorphism. Since being an epimorphism is a local property, without loss of generality we can assume that $X$ is affine and the line bundles $N_{X/Y}$ and $N_\F$ are trivial. By definition, a section $h\in \I_{\Per}$ must verify $[hd\omega] = 0\in \Omega^2_{X}/(\omega\wedge\Omega^1_X)$. Equivalently, there is some $\theta\in \Omega^1_X$ such that $hd\omega=\omega\wedge\theta$. Setting $\eta=\theta+dh$, we get $hd\omega - \omega\wedge(\eta - dh) = 0$, and then $(h,[\eta])\in\PP_{\F,N_{X/Y}^{\otimes - n}}^{\,\Int}$.
\end{proof}

In particular, combining this last proposition with Proposition \ref{obst:h1_pint} we see that the obstructions considered previously actually lie in $H^1(X,\I_{\Per} \otimes N_\F\otimes N_{X/Y}^{\otimes - n})$. We will end this section by stating a criterion for the vanishing of these groups, taking into consideration the positivity of the embedding and the geometry of $\Per(\F)$. For a criterion regarding only the former, the reader is referred to Section \ref{section:positivity}.

\begin{proposition}
Let $X\subseteq Y$ be a smooth compact divisor and $\F$ a foliation on $X$ such that $\Per(\F)$ is reduced and without isolated points and $N_\F\otimes N_{X/Y}^{\otimes - n}$ is anti-ample for every $n\geqslant 1$. Then $H^1\big(X,\I_{\Per}\otimes N_\F\otimes N_{X/Y}^{\otimes -n}\big)=0$ for every $n\geqslant 1$.  
\end{proposition}
\begin{proof}
The anti-ampleness of the line bundles implies that $H^1(X,N_\F\otimes N_{X/Y}^{\otimes - n})=0$. Moreover, since $\Per(\F)$ does not have isolated points this also implies that $h^0(\Per(\F),N_\F\otimes N_{X/Y}^{\otimes - n})=0$. Looking at the cohomology of the exact sequence
$$ 0\rightarrow \I_{\Per}\otimes N_\F\otimes N_{X/Y}^{\otimes - n}\rightarrow N_\F\otimes N_{X/Y}^{\otimes - n}\rightarrow  \OO_{\Per}\otimes N_\F\otimes N_{X/Y}^{\otimes - n}\rightarrow 0$$
we get the desired vanishing.
\end{proof}

\subsection{Unobstructedness of foliation singularities} 
In this section, we will analyze the (un)obstructedness of different types of singularities. In particular, we will show a large class of unobstructed singularities, justifying the assumptions made in the previous section.

\begin{example}{\label{unf_int1}}
Let $\F$ be a foliation with a local holomorphic first integral $f$ at $0\in \C^m$. Writing $\omega=df$, we can describe the extension space as
\[
\Ex(\F) \simeq \frac{\OO_{\C^m,0}}{(\partial f/\partial z_1,\dots,\partial f/\partial z_m)}.
\]
Observe that if $0$ was a regular point of $\F$ then we automatically get $\Ex(\F)=0$. Since
$$ d(f+t h)= df+t dh + hdt,$$
we see that every first order unfolding of $\F$ is of this form. 
Moreover, we claim that any infinitesimal unfolding of $\F$ is of the form $\omega^{(n)}=d(f_n)$ for some element $f_n\in \OO_{X^{(n)}}$ extending $f$.
Indeed, let $\omega^{(n)}$ be an extension of $\omega_\F$. Proceeding by induction, we can freely assume that the restriction $\omega^{(n-1)}=\omega^{(n)}\vert_{X^{(n-1)}}$ is of the form $\omega^{(n-1)}=df_{n-1}$. In particular, we can construct an extension $\overline{\omega}^{(n)}=df_n$. By Proposition \ref{def:fol}, we see that the isomorphism class of $\omega^{(n)}$ relative to $X^{(n-1)}$ must be the same as the class of
$$\overline{\omega}^{(n)}+ d(t h)$$
 and therefore $\omega^{(n)}$ is of the desired type. In particular, foliations' singularities with first integrals are unobstructed. By a theorem of Malgrange, this includes singularities with $\codim(\Sing(\F))\geq 3$.
\end{example}

A particular case of unobstructedness is the following:

\begin{definition}
    A germ of foliation singularity $\F$ is said to be \textit{rigid} if $\Ex(\F)=0$.
\end{definition}

\begin{lemma}
    Let $\F$ be a \textit{rigid} germ of foliation singularity. Then every infinitesimal extension of $\F$ is isomorphic to the trivial extension. 
\end{lemma}
\begin{proof}
In view of Proposition \ref{def:fol}, this follows easily by induction on the order of the extension.
\end{proof}

\begin{example}{\label{unf_kupka}}
Let $\F$ be a Kupka singularity on $(\C^n,0)$. By Kupka's Theorem \cite[pag.197]{neto2020complex} there must be a coordinate system $z_1,\dots,z_n$ around the origin and a representative of $\F$ of the form $\omega = a_1dz_1+ a_2dz_2$ with $a_1$ and $a_2$ depending only on the first two variables. Let $h$ be an element in $\I_{\Per}$ and $\sigma$ a local $1$-form such that $hd\omega = \omega\wedge\sigma$. Observe that the coefficients of $dz_k$ in $\sigma$ must vanish for $3\leqslant k \leqslant n$, i.e., $\sigma = b_1 dz_1 + b_2 dz_2$. Since $d\omega$ does not vanish at the origin,
\begin{align*}
 hd\omega = \omega\wedge\sigma \quad &\Rightarrow \quad h(\partial a_2/\partial z_1 - \partial a_1/\partial z_2) dz_1\wedge dz_2 = (a_2b_1-a_1b_2)dz_1\wedge dz_2\\
 & \Rightarrow \quad h(\partial a_2/\partial z_1 - \partial a_1/\partial z_2) = (a_2b_1-a_1b_2)\\
 & \Rightarrow \quad h = (\partial a_2/\partial z_1 - \partial a_1/\partial z_2)^{-1}(a_2b_1-a_1b_2).
\end{align*}
This is, $h$ is the contraction of $\omega$ with $(\partial a_2/\partial z_1 - \partial a_1/\partial z_2)^{-1}(b_1\partial z_2-b_2\partial z_1)$. This shows that Kupka singularities are rigid, hence unobstructed.
\end{example}

\begin{lemma} \label{ex:locallyfree}
Suppose $\F$ is a codimension $1$ foliation on $(\C^n,0)$ with locally free tangent sheaf whose non-Kupka singularities have codimension at least 3. Then $\Ex(\F) = 0$.
\end{lemma}

\begin{proof}
We follow \cite[Theorem 3.1]{cerveau1991determinacy} closely. Let $U\subseteq \C^n$ be an open Stein neighborhood of the origin and $\omega\in\Gamma(\Omega^1_U)$ a representative of $\F$. Take $\U = \{U_i\}_i$ a Stein open cover of $W= U\setminus\Sing(\omega)\cap\Sing(d\omega)$. Taking Examples \ref{unf_int1} and \ref{unf_kupka} into account, if $hd\omega$ is a multiple of $\omega$ we can find vector fields $X_i\in \Gamma(T_{U_i})$ such that $h = i_{X_i}\omega$. In particular $i_{X_i-X_j}\omega = 0$, and therefore $\{X_i-X_j\}_{ij}$ is a $\smash{\check{C}}$ech $1$-cocycle with coefficients in $T_\F$. Standard vanishing theorems imply $H^1(W,T_\F) = 0$ \cite[pág. 133]{stein}, so we can obtain local sections $Y_i\in \Gamma(T_\F|_{U_i})$ verifying $X_i-X_j= Y_i - Y_j$. Therefore, there is a vector field $X$ over $W$ tangent to the foliation such that $X|_{U_i} = X_i-Y_i$. By Hartogs' Theorem it can be extended to a field $X\in\Gamma(T_U)$ with $h = i_X\omega$, finishing the proof.
\end{proof}

Rigidity can also be interpreted in terms of persistent singularities. The following lemma is illustrative and encompasses some of the examples above.

\begin{lemma} \label{lemma:per=rigid}
Let $\F$ be a germ of foliation. Then $\F$ is rigid if and only if $\Per(\F)=\Sing(\F)$.
\end{lemma}
\begin{proof}
It suffices to observe that the singular ideal of $\F$ is generated by the contractions $\imath_v(\omega_\F)$ with $v\in T_X$, while
$ I_{\Per}(\F)=\{ h \, \vert \, \omega_\F+ t\eta +hdt \mbox{ is integrable on } X\times D \}$.
Then the lemma follows from Remark \ref{trivialextensions} and the uniqueness of $[\eta]$ established in Lemma \ref{eq:integrabilidad}.
\end{proof}

\section{Formal and algebraic extensions}
\subsection{Algebraization of infinitesimal extensions}\label{section4.1}

We will now proceed to prove the main results as stated in the introduction. The first step is to construct formal extensions, as defined in Definition \ref{def:formalext}. The following theorem is a straightforward application of Proposition \ref{obst:h1_pint} and Proposition \ref{pint=per}.

\begin{theorem} \label{ext:formal}
Let $X\subseteq Y$ be a regularly embedded divisor. If $\F$ is a foliation on $X$ having unobstructed singularities and such that 
$$H^1\big(X,\I_{\Per}\otimes N_\F\otimes N_{X/Y}^{\otimes -n}\big)=0$$ 
for every $n \geqslant 1$, then there is a formal extension $\widehat{\F}$ of $\F$.
\end{theorem}

Having already settled the existence of formal extensions, we will now turn our attention to their algebraicity. First, let us briefly recall the Lefschetz conditions. For further details the reader is referred to \cite[Chapter 10]{badescu2012projective}.

\begin{definition} \label{lefschetz}
   Let $X\subseteq Y$ be a closed subvariety. We say that $(X,Y)$ satisfies the \textit{Lefschetz condition} $\Lef(X,Y)$ if for every open set $U$ containing $X$ and every locally free sheaf $E$ on $U$, the natural map 
   $$H^0(U,E)\to H^0\big(\,\widehat{X}_{/Y},i^*E\big)$$
   is an isomorphism. The pair $(X,Y)$ satisfies the \textit{effective Lefschetz condition} $\Leff(X,Y)$ if $\Lef(X,Y)$ holds and every locally free sheaf $\E$ on $\widehat{X}_{/Y}$ is of the form $i^*E$ for some locally free sheaf $E$ on some neighborhood $U$ of $X$.
\end{definition}

\begin{theorem}\label{ext:global}
   Let $\F$ be a foliation on a smooth variety $X$ having unobstructed singularities and let $X\subseteq Y$ be a regular embedding with $\dim(Y)=\dim(X)+1$.  
   Suppose that $$H^1(X,\I_{\Per}\otimes N_\F\otimes N_{X/Y}^{\otimes -n})=0$$
   for every $n\geqslant 1$ and $\Leff(X,Y)$ holds.
    Then $\F$ extends to a foliation on $Y$. 
\end{theorem}

\begin{proof}
Our hypotheses imply that the problem of extending $\F$ to the infinitesimal neighborhoods $X^{(k)}$ is unobstructed. This is, we can construct line bundles $\L^{(k)}$ and global sections $\omega^{(k)}$ of $\Omega^1_{X^{(k)}}\otimes \L^{(k)}$ compatible with their respective pullbacks. 
The property $\Leff(X,Y)$ guarantees the existence of a (Zariski) neighborhood $U$ of $X$ and a line bundle $\overline{\L}$ on $U$ such that $\overline{\L}\vert_{X^{(k)}}\simeq \L^{(k)}$ in a compatible manner. For the sake of clarity, we will regard all these sheaves as $\OO_U$-modules via the corresponding pushforwards, which we will omit in order to keep the notation manageable.
Let us first show that sheaves $\Omega^1_{X^{(k)}}\otimes \L^{(k)}$ form an inverse limit such that the pullback of twisted differential forms induces an isomorphism
\begin{equation}
    \label{isolim} \lim_{\longleftarrow}\Omega^1_U\otimes \overline{\L} \otimes \OO_{X^{(k)}} \simeq \lim_{\longleftarrow} \Omega^1_{X^{(k)}} \otimes \L^{(k)}.
\end{equation}
Observe that for every $k$ we have the exact sequence
$$C^{(k)}: 0\rightarrow \mathcal G^{(k)}\rightarrow \Omega^1_U\otimes\overline{\L}  \otimes \OO_{X^{(k)}} \rightarrow \Omega^1_{X^{(k)}} \otimes \L^{(k)} \rightarrow 0,$$
where $\mathcal G^{(k)}\subseteq \Omega^1_U\otimes\overline{\L}$  corresponds to the subsheaf generated by elements of the form $(dt^{k+1})$, for $t$ a local equation for $X$ in $Y$.
For $k<k'$ we have (compatible) maps of exact sequences of $\OO_U$-modules
$$C^{(k')}\to C^{(k')}\otimes \OO_{X^{(k)}}\to C^{(k)}$$ 
of the form
\begin{center}
\begin{tikzcd}
\mathcal G^{(k')} \arrow[r] \arrow[d] & \Omega^1_U \otimes \overline{\L} \otimes \OO_{X^{(k')}}\arrow[r] \arrow[d] & \Omega^1_{X^{(k')}}\otimes \L^{(k')}\arrow[d] \\
    \mathcal G^{(k')} \otimes \OO_{X^{(k)}} \arrow[r] \arrow[d,"0"]& \Omega^1_U \otimes \overline{\L} \otimes \OO_{X^{(k)}}\otimes \OO_{X^{(k')}}\arrow[r] \arrow[d] &  \Omega^1_{X^(k')}\otimes \L^{(k')}\otimes \OO_{X^{(k)}} \arrow[d,"i^*"] \\
    \mathcal G^{(k)} \arrow[r] & \Omega^1_U \otimes \overline{\L} \otimes \OO_{X^{(k)}}\arrow[r] & \Omega^1_{X^{(k)}}\otimes \L^{(k)}.
\end{tikzcd}
\end{center}
We will denote by $z_{kk'}$ (resp. $r_{kk'}$, $i^*_{kk'}$) the morphism given by the composition of the arrows in the first (resp. second, third) column.
All this adds up to an exact sequence of inverse systems of sheaves of $\OO_{U}$-modules
$$ 0\rightarrow \{( \mathcal G^{(k)}),z_{kk'}\} \rightarrow \{\Omega^1_U\otimes \overline{\L} \otimes \OO_{X^{(k)}},r_{kk'}\} \rightarrow  \{ \Omega^1_{X^{(k)}} \otimes \L^{(k)},i^*_{kk'}\} \rightarrow 0 .$$
Now since $z_{kk'}=0$ for every $k$ and $k'$, the leftmost inverse system satisfies both $$\lim_{\longleftarrow} \mathcal G^{(k)} =0$$ and the Mittag-Leffler condition (see \cite[Section 0.13.2]{EGAIIIp1}).
In particular, by \cite[Proposition 0.13.2.2]{EGAIIIp1} we get the desired isomorphism (\ref{isolim}).
With this in mind, observe that the twisted differential forms $\omega^{(k)}$ define a global section of $\varprojlim \Omega^1_{X^{(k)}} \otimes \L^{(k)}$ and therefore can be understood as a compatible system of morphisms 
$$ \OO_{X^{(k)}} \rightarrow \Omega^1_U\otimes \overline{\L}\otimes \OO_{X^{(k)}}.$$
Again, by $\Leff(X,Y)$ the formal family of global sections must be induced by some 
$$\overline{\omega}\in H^0(U,\Omega^1_U\otimes \overline{\L}).$$ 
The integrability of $\overline{\omega}$ now also follows from $\Leff(X,Y)$, since what we have done so far shows that $\smash{\overline{\omega}\wedge d \overline{\omega}: \OO_U \to \Omega^3_U\otimes \mybar{1}{0pt}{\L}^{\otimes 2}}$ restricts to zero on each $X^{(k)}$.
This is, this twisted differential form defines a foliation on $Y$ that (by construction) extends $\F$. 
\end{proof}

As a first application, let us consider the setting where $Y$ is a smooth projective variety.
By a theorem of Grothendieck (\cite[Theorem 10.2]{badescu2012projective}), every ample divisor $X\subseteq Y$ satisfies $\Leff(X,Y)$. We automatically get the following. 

\begin{corollary} \label{ext_ample}
    Let $X$ be a smooth ample divisor on a smooth projective variety $Y$ and $\F$ a codimension $1$ foliation on $X$ with unobstructed singularities and $$H^1(X,\I_{\Per}\otimes N_\F\otimes N_{X/Y}^{\otimes -n})=0$$ for every $n\geqslant 1$. Then $\F$ extends to a foliation on $Y$.
\end{corollary}

\begin{example}
    Let $X\subseteq Y$ be a smooth ample divisor with $\dim(X)\geqslant 3$, and $\L_1,\dots, \L_k$ ample line bundles on $X$. A \emph{logarithmic foliation} on $X$ is given by a differential form
    $$ \omega_\F= f_1\dots f_k\sum_{i=1}^k \lambda_i \frac{df_i}{f_i},$$
    for some sections $f_i\in H^0(X,\L_i)$, and complex numbers $\lambda_i\in\C\setminus\{0\}$ verifying the condition \smash{$\sum_{i=1}^k \lambda_k c_1(\L_k) = 0$}. This last equation is needed to assure that $\omega_\F$ defines a global twisted differential form \cite[p. 752]{calvo1994irreducible}. Following \cite{cukierman2006singularities}, if the sections $f_i$ are generic, the singular scheme of $\F$ is  of the form 
    $$\Sing(\F)=\overline{K(\F)} \cup R,$$
    where $R$ is a collection of isolated points and $\overline{K(\F)}=\cup_{i\neq j} V(f_i,f_j)$.
    The local computation carried out in \cite[p. 134]{cukierman2006singularities} implies that $T_\F$ is locally free around every point of $\overline{K(\F)}$, and therefore by Lemma \ref{ex:locallyfree} we can assert that those singularities are rigid (hence unobstructed). The remaining isolated singularities in $R$ are also unobstructed by Example \ref{unf_int1}. Moreover, since in this case we have $\Per(\F)=\overline{K(\F)}$ the ideal of persistent singularities is generated by the elements $\widehat{f}_i=f_1\dots f_{i-1}f_{i+1}\dots f_k$. In particular, setting $\L=(\L_1\otimes \cdots \otimes \L_k)$ and $\widehat{\L}_i=\L\otimes \L_i^{-1}$ we have a resolution
    $$ 0 \rightarrow (\L^{-1})^{\oplus (k-1)} \longrightarrow \bigoplus_{i=1}^k \widehat{\L}_i^{ -1} \xrightarrow{\widehat{f}_1,\dots ,\widehat{f}_k} \I_{\Per} \rightarrow 0,  $$
    where the first map is given by the trivial relations $f_i\widehat{f}_i-f_j\widehat{f}_j=0$. Observe that $N_\F\simeq \L$, and hence tensoring by $N_\F\otimes N_{X/Y}^{\otimes -n}$ and looking at the cohomology sequence we get 
    $$H^1(X,\I_{\Per}\otimes N_\F\otimes N_{X/Y}^{\otimes -n} )\simeq \bigoplus_{i=1}^k H^1(X,\L_i\otimes N_{X/Y}^{\otimes -n}).$$
    This isomorphism indicates that the obstructions to extend the foliation are precisely the obstructions to extend the sections $f_i$ \cite[Prop. 3.3.4]{sernesi2007deformations}. In particular, this implies that $\F$ can be extended to $Y$ whenever $\L_i\otimes N_{X/Y}^{\otimes -1}$ is anti-ample or $X$ is arithmetically Cohen-Macaulay.
\end{example}

\begin{corollary} \label{cor:unfoldings}
    Let $\F$ be a codimension $1$ foliation on a smooth proper variety $X$ having unobstructed singularities and satisfying $$H^1(X,\I_{\Per}\otimes N_\F)=0.$$ Then $\F$ admits an unfolding along any deformation $\X\to \Spec(\C[[t]])$ of $X$. 
\end{corollary}
\begin{proof}
    The statement follows directly from Theorem \ref{ext:global} and the fact that Grothendieck's Extension Theorem (\cite[Theorem 5.1.4]{EGAIIIp1}) -- or \cite[Theorem 1.2]{kosarew1991grothendieck} in the analytic case -- ensures $\Leff(X,\X)$.
\end{proof}

\begin{example}
Let $X\subseteq \P^{m}$ be a smooth hypersurface. If $m\geq 3$, computation shows that every deformation $\X\to \Spec(\C[[t]])$ of $X$ is an embedded deformation. Although this is probably the easiest example of a family of varieties over $T$, the existence and characterization of unfoldings along these deformations is non-trivial. Of course, it may happen that the foliation $\F$ on $X$ extends to a foliation $\overline{\F}$ on $\P^{m}$ (see Corollary \ref{ext_ample} or Section \ref{section:positivity} below).  In this case, the restriction of the trivial unfolding of $\overline{\F}$ along $\P^{m+1} \times \Spec(\C[[t]])$ will give an unfolding along any (embedded) deformation of $X$. If $\F$ does not extend, we can still aim to construct unfoldings along deformations. For instance, if the singular scheme of $\F$ is reduced (implying $\Per(\F)=\overline{K(\F)}$) and the non-Kupka components of $\Sing(\F)$ are of codimension at least three and do not intersect $\overline{K(\F)}$, then $\F$ will have unobstructed singularities. If moreover $\Per(\F)$ satisfies $H^1(\I_{\Per(\F)/\P^{m}}\otimes N_\F)=0$, then tensoring the exact sequence
$$ 0\rightarrow I_{X/\P^{m}}\rightarrow \I_{\Per(\F)/\P^{m}}\rightarrow i_*\I_{\Per(\F)}\rightarrow 0
$$
with $N_\F$ and taking cohomology on $\P^{m}$ we get  $H^1(X,\I_{\Per}\otimes N_\F)=0$ and therefore we are still able to unfold $\F$ along any deformation of $X$. 
\end{example}

\begin{remark}
     Under the hypothesis of Corollary \ref{cor:unfoldings}, one can further show that  $\F$ admits a non-trivial unfolding along $X\times \Spec(\C[[t]])\to \Spec(\C[[t]])$ whenever $\Ex(\F)\neq 0$. Indeed, by the unobstructedness of $\F$ one can construct a formal family which is non-trivial at first order. By the argument above, this family can be algebraized to a foliation $\overline{\F}$ on $X\times \Spec(\C[[t]])$. The triviality of $\overline{\F}$ would contradict \cite[Theorem 5.1.4]{EGAIIIp1}. 
\end{remark}

\subsection{A positivity criterion for extendability}\label{section4.2}
\label{section:positivity}
A large difference in positivity between the embedding $X\hookrightarrow Y$ and the normal bundle of the foliation can provide sufficient conditions for extendability.
This is the spirit of \cite{figueira2023extensions}, where the author shows that a foliation $\F$ on a smooth hypersurface $X\subseteq \P^{n+1}$ with $\dim X\geqslant 3$ extends (uniquely) provided that $\deg(X)> 2\deg(\F)+1$. 

Recall from Definition \ref{def:Iper} that the sheaf of ideals $\I_{\Per}$ is defined as the kernel of the multiplication map \smash{$d\omega_\F:\OO_X\rightarrow \Omega_X^2(N_\F)/(\omega_\F\wedge\Omega^1_X)$}. In particular, the twisted ideal of persistent singularities appearing in the Corollary \ref{ext_ample} fits in the exact sequence
\[
0\longrightarrow \I_{\Per}\otimes N_\F\otimes N_{X/Y}^{\otimes - n}\longrightarrow N_\F\otimes N_{X/Y}^{\otimes - n}\longrightarrow  Q_{\F}(n)
\]
where the sheaf on the right was defined in (\ref{def:QF}) as the quotient
\[
0\longrightarrow \Omega^1_X(N_\F\otimes N_{X/Y}^{\otimes - n})\xrightarrow{\omega_\F\wedge}\Omega^2_X(N_\F^2\otimes N_{X/Y}^{\otimes - n})\longrightarrow  Q_{\F}(n)\longrightarrow 0.
\]
The hypothesis $H^1(X,\I_{\Per}\otimes N_\F\otimes N_{X/Y}^{\otimes -n})=0$ then follows from the conditions
\begin{equation}\label{eq:2}
H^1\big(X,N_\F\otimes N_{X/Y}^{\otimes -n}\big) = H^1\big(X,\Omega^1_X(N_\F\otimes N_{X/Y}^{\otimes -n})\big) = H^0\big(X,\Omega^2_X(N_\F^{\otimes 2}\otimes N_{X/Y}^{\otimes -n})\big) = 0.   
\end{equation}
Let us now suppose that \smash{$N_{X/Y}$} is ample. In this situation, if \smash{$N_\F^{\otimes 2}\otimes N_{X/Y}^{\otimes -n}$} is anti-ample, then so is \smash{$N_\F\otimes N_{X/Y}^{\otimes -n}$}, since 
$$(N_\F\otimes N_{X/Y}^{\otimes -n})^{\otimes 2}\simeq (N_\F^{\otimes 2}\otimes N_{X/Y}^{\otimes -n})\otimes N_{X/Y}^{\otimes -n}.$$ Applying Kodaira-Nakano Vanishing Theorem to both $N_\F^{\otimes 2}\otimes N_{X/Y}^{\otimes -n}$ and $N_\F\otimes N_{X/Y}^{\otimes -n}$ we get:

\begin{theorem} \label{teo:positivity criterion}
Let $X$ be a smooth ample divisor on a smooth projective variety $Y$ and $\F$ a codimension $1$ foliation on $X$ with unobstructed singularities. If $\dim X\geqslant 3$ and $N_\F^{\otimes 2}\otimes N_{X/Y}^\vee$ is anti-ample, then $\F$ extends to a foliation on $Y$.
\end{theorem}

\begin{example} \label{ej:figueira}
Suppose $X\subseteq \P^{m+1}$ is a smooth hypersurface with $m\geqslant 3$ and $\deg(X)\geqslant 3$. If $\F$ is a codimension $1$ foliation on $X$ satisfying the inequality $$\deg(X)> 2\deg(\F)+1,$$ then the cohomology groups in equation (\ref{eq:2}) vanish \cite[Lemma 4.5]{araujo2013fano}. Therefore, if the singularities of the foliation are unobstructed, $\F$ extends to a foliation on $\P^{m+1}$. Then we recover the main theorem in \cite{figueira2023extensions} if this additional hypothesis is assumed. Moreover, this extension is in fact unique (see Proposition \ref{uniqueness} below).
\end{example}

\section{The foliated Zak-L'vovsky problem}\label{section5}

To continue with the study of extendability of foliations we will make a brief digression into the extendability problem of subvarieties of projective spaces. Let $X\subseteq \P^m$ be a closed subvariety. Given a linear embedding $\P^m\subseteq \P^{m+1}$, an \textit{extension} of $X$ is a closed subvariety $Y\subseteq \P^{m+1}$ of dimension $\dim Y = \dim X + 1$ such that $X = Y \cap \P^m$, where this intersection is taken in the scheme-theoretical sense. If we choose a point $p\in\P^{m+1}\setminus\P^m$, the projective \textit{cone} of $X$ with vertex $p$ is clearly an extension. It is natural to ask whether these are the only possible extensions. This question was posed by Zak and L'vovsky and they arrived to the following conclusion:

\begin{theorem}{\cite[Coro. 1]{l1992extensions}}\label{L'vovsky}
    Let $X\subseteq \P^m$ be a smooth subvariety of dimension $\dim(X)\geqslant 2$ that is not a quadric. If $\smash{H^1\big(X,T_X(-1)\big) = 0}$ then every extension of $X$ to $\P^{m+1}$ is a cone.
\end{theorem}

Those definitions have a natural analogy in the foliated setting. Given a linear embedding $\P^m\subseteq \P^{m+1}$, there is a trivial way to extend a foliation $\F$ on $\P^n$ to $\P^{m+1}$: just pick a linear retraction $\pi:\P^{m+1}\dashrightarrow \P^m$ and consider the pullback $\pi^\ast\F$. This is the foliation whose leaves are cones with some vertex $p\in \P^{m+1}\setminus \P^m$ of the leaves of $\F$. Such foliation is called a \textit{cone} of $\F$.\\ 

It is reasonable to seek a foliated version of Theorem \ref{L'vovsky}. To do so, we need a replacement for the vanishing condition on $H^1(X,T_X(-1))$. This last group parametrizes a special class of infinitesimal extensions of $X$ called \textit{thickenings}, bent according to normal bundle $N_{X/\P^{m}}=\OO_X(1)$ ( see \cite[p. 15-16]{sernesi2007deformations}). Fortunately, the analogy can be pushed further utilizing the tools at our disposal.

\begin{definition} 
Suppose that $\F$ is a codimension $1$ foliation on $\P^m$ with normal bundle $\OO_{\P^m}(k)$. We say that $\F$ is \textit{Camacho-Lins Neto regular} if
\[
\Ex^n(\F,\OO_{\P^m}(1)) = 0 \quad\forall\; n < k. 
\]
\end{definition}

At first sight it is not clear that this notion is the same as the original definition stated in \cite{camacho1982topology}, but they turn out to be equivalent \cite[pag. 1600/08]{molinuevo2016unfoldings}. The next result is implicit in \cite[Prop. 3]{camacho1982topology}

\begin{theorem} \label{teo:CLN regularity}
    Let $\F$ be a codimension $1$ foliation on $\P^m$. The following are equivalent:
\begin{enumerate}
    \item $\F$ is Camacho-Lins Neto regular.
    \item given a linear embedding $\P^m\subseteq \P^{m+1}$, every extension of $\F$ to $\P^{m+1}$ is a cone.
\end{enumerate}
\end{theorem}

\begin{proof}
If $N_\F=\OO_{\P^m}(k)$, an extension of $\F$ to $\P^{m+1}$ is a cone if and only if its restriction to every infinitesimal neighborhood $X^{(k)}$ of $X=\P^m\subseteq \P^{m+1}$ is isomorphic to the trivial extension. But this holds if and only if the restrictions of such an extension of $\F$ to each intermediate neighborhood $$\P^m = X \subseteq X^{(1)}\subseteq X^{(2)}\subseteq \dots \subseteq X^{(k-1)}\subseteq X^{(k)}\subseteq \P^{m+1}$$ are trivially extended at each step. Since at each step these extensions are classified by the groups $\Ex^n(\F,\OO_{\P^m}(1))$ with $n< k$, respectively, the theorem follows.
\end{proof}

\begin{example}
Suppose that $\F$ is a foliation on $\P^n$ with split tangent sheaf and whose non-Kupka singularities have codimension at least 3. One of the main results of \cite{massri2018kupka} implies
\[
\Ex^n(\F,\OO_{\P^m}(1))=0 \quad\forall\, n,
\]
and therefore $\F$ is Camacho-Lins Neto regular. This provides a family of foliations whose only extensions are cones.
\end{example}

The work done in the previous sections allows us to prove a similar result for ambient spaces other than $\P^m$.

\begin{proposition}\label{uniqueness}
    Let $X$ be a smooth ample divisor on a smooth projective variety $Y$. Assume $X\subseteq Y$ is not the inclusion of a hyperplane in a projective space, and $\F$ is a codimension $1$ foliation on $X$ such that
    $$ \Ex^n(\F,N_{X/Y})=0$$
    for every $1\leqslant n< \min\big\{ d\geqslant 1 \,\vert \, N_\F\otimes N_{X/Y}^{\otimes -d} \mbox{ is anti-ample}\big\}$ . If $\F_1$ and $\F_2$ are two extensions of $\F$ to $Y$, then $\F_1 = \F_2$.
\end{proposition}
\begin{proof}
Since the vanishing for the higher extension spaces follows from the anti-ampleness of $N_\F\otimes N_{X/Y}^{\otimes -n}$ we can conclude that $\Ex^n(\F,N_{X/Y})=0$ for every $n$. By the Mori-Sumihiro-Wahl Theorem \cite{Wahl} we have $H^0(X,T_X\otimes N_{X/Y}^{\otimes -n})=0$. This implies that if $\F_1$ and $\F_2$ are two extensions of $\F$, then one can easily see (by induction on $n$) that their restriction to $X^{(n)}$ actually coincide. But then $\Leff(X,Y)$ implies that $\F_1=\F_2$.
\end{proof}

\begin{remark}
    This last argument also indicates how to construct a wide variety of foliations that do not extend. Indeed, it may happen that a non-integrable codimension $1$ distribution $\D$ on $Y$ restricts to a foliation $\F$ on $X$, where $X\subseteq Y$ is an ample divisor (this happens trivially, for instance, if $X$ is a surface). Let us further suppose that $N_\F\otimes N_{X/Y}^{\otimes -1}$ is anti-ample. In this case, easy computation shows that \smash{$H^0(X,\PP_{\F,N_{X/Y}^{\otimes - n}})=0$}. By Proposition \ref{extensions:dist} and the argument above, this implies that $\F$ has a unique extension to $Y$ as a distribution. In particular, it does not extend to a foliation on $Y$.
\end{remark}

\section{Rigidity theorems}\label{section6}
As mentioned in the introduction, one reason a foliation might extend is the existence of tubular neighborhoods.

\begin{definition}
 Let $X\subseteq Y$ be a regular embedding of codimension $1$. A \emph{tubular neighborhood} of $X$ in $Y$ is an analytic open subset $U\subseteq Y$ containing $X$ and admitting a retraction $r:U\to X$. 
 A \emph{Zariski tubular neighborhood} is just a tubular neighborhood where $U$ is a Zariski open and the retraction $r$ is algebraic.
 We say that the $n$-th infinitesimal neighborhood of $X$ in $Y$ is tubular if there exists a retraction $X^{(n)}\to X$. 
\end{definition}

 The existence of such neighborhoods is, however, a non-trivial problem. Griffiths showed in \cite[Proposition 1.6]{griffiths1966extension} that the $n$-th obstructions to finding tubular neighborhoods of $X$ in $Y$ lies in the cohomology group
\[
H^1\big(X,T_X\otimes N_{X/Y}^{\otimes -n}\big).
\]
This is intuitively reasonable because these cohomology groups are related to the splitting of the normal sequence of $X\subseteq Y$. Several theorems in the literature provide tubular neighborhoods assuming the vanishing of those groups and some additional ampleness hypothesis (see for instance \cite[Corollary p.363]{grauert1962modifikationen} and \cite{fujita1983rational}). This, together with what we have done so far in this work, provides us with enough tools to give a foliated version of these statements.

\begin{theorem} \label{fol_grauert}
    Let $X$ be a smooth projective variety with $\dim(X)\geq 3$ and $X\subseteq Y$ a regular embedding with ample normal bundle. Suppose further that there exists a codimension $1$ foliation $\F$ on $X$ with rigid singularities admitting an extension to $Y$ and satisfying
    $$H^1\big(X,T_\F\otimes N_{X/Y}^{\otimes -n}\big) = 0$$
for each $n\geqslant 1$. Then $X$ has a tubular neighborhood. If moreover $Y$ is projective and $X$ is an ample divisor, then $X$ admits a Zariski tubular neighborhood and we have the following dichotomy:
\begin{enumerate}
    \item the extension of $\F$ is unique, or
    \item $X\subseteq Y$ is a hyperplane included in a projective space and the extension of $\F$ is a cone.
\end{enumerate}
\end{theorem}

\begin{proof} 
By Lemma \ref{lemma:per=rigid} we get $\Sing(\F) = \Per(\F)$.
Tensoring the exact sequence
\[
0\rightarrow T_\F \rightarrow T_X \rightarrow \I_{\Per} \otimes N_\F\rightarrow 0
\]
by $N_{X/Y}^{\otimes -n}$ and looking at its cohomology we see that
\begin{enumerate}
    \item $H^0\big(X,T_X\otimes N_{X/Y}^{\otimes -n}\big)\rightarrow H^0\big(X,\I_{\Per}\otimes N_\F\otimes N_{X/Y}^{\otimes -n}\big)$ is surjective, and
    \item $H^1\big(X,T_X\otimes N_{X/Y}^{\otimes -n}\big)\rightarrow H^1\big(X,\I_{\Per}\otimes N_\F\otimes N_{X/Y}^{\otimes -n}\big)$ is injective.
\end{enumerate}
Moreover, the arrow in $(2)$ maps the $n$-th obstruction to find a tubular neighborhood to the obstructions to extend the foliation $\F$. Note that the vanishing of the former naturally implies the vanishing of the latter, because $\F$ extends as the pullback by the retraction of infinitesimal tubular neighborhoods. Since $\F$ extends to $Y$, these elements do vanish. By \cite[Theorem II]{griffiths1966extension} and the following remark we can conclude that $X$ has a tubular neighborhood. 

If moreover $Y$ is projective and $X\subseteq Y$ is an ample divisor, then - just as in Corollary \ref{ext_ample} - we know that $\Leff(X,Y)$ holds. 
Let us fix a projective embedding $X\subseteq \P^N$. With this in mind, the retraction $r$ defining the tubular neighborhood can be written in the form  $r=[s_0:\cdots:s_N]$ where $s_i$ is a section of the line bundle $\L:=r^*\OO_X(1)$. Since the restrictions of the sections $s_i$ to the infinitesimal neighborhoods of $X$ define a section of the formal completion of the line bundle $\L$ along $X$, the effective Lefschetz condition  $\Leff(X,Y)$ allows us to extend both $\L$ and the sections $s_i$ to a Zariski neighborhood of $X$. These new sections restricted to the complement of their common zeroes provide a Zariski tubular neighborhood.
On the other hand, (1) implies that $$\Ex^n(\F,N_{X/Y})=0$$ for every $n\geq 1$. The conclusion now follows from Theorem \ref{teo:CLN regularity} and Proposition \ref{uniqueness}.
\end{proof}

A hypersurface of $\P^m$ has a tubular neighborhood if and only if it is a hyperplane, by Van de Ven's Theorem \cite[Theorem 7.1]{badescu2012projective}. We can therefore arrive at the following conclusion:

\begin{corollary}
Suppose $X\subseteq \P^{m}$ is a smooth hypersurface of dimension at least 3 and $\deg(X)\geqslant 2$. If $\F$ is a codimension one foliation with rigid singularities and $H^1\big(X,T_\F\otimes N_{X/Y}^{\otimes -n}\big) = 0$ for every $n\geqslant 1$, then $\F$ does not extend.
\end{corollary}

Moving to the unfolding scenario, \cite[Theorem 5.5]{gomez1989unfoldings} states that every unfolding of a foliation with $H^1(X,T_\F)= 0$ and rigid reduced singularities on a smooth surface is analytically trivial. With our approach, we can easily generalize this to the following.

\begin{theorem} \label{unf:trivial}
Suppose $\F$ is a codimension $1$ foliation with rigid singularities on a proper smooth variety $X$. If $H^1(X,T_\F) = 0$, then every unfolding of $\F$ is analytically trivial. 
\end{theorem}

\begin{proof}
Given a deformation $X^{(n-1)}$ of order $(n-1)$ of $X$, the space $H^1(X,T_X)$ acts freely and transitively on the space of deformations of order $n$ extending $X^{(n-1)}$ (see for instance \cite[Proposition 1.2.12]{sernesi2007deformations}). Proceeding as in the proof of Theorem \ref{fol_grauert} we see that such unfolding does not deform $X$ and $\Ex(\F) = 0$, hence every formal unfolding over $\C[[t]]$ of $\F$ is trivial. 

For the remainder of the proof we will follow the spirit of  \cite[page 514]{gomez1989unfoldings}.
Let $\dF$ be an unfolding of $\F$ along a deformation $\X\to \Delta$ of $X$ over a small disc $\Delta\subseteq\C$. By the argument above, $\pi:\X\to \Delta$ is formally trivial and hence analytically trivial by \cite{schuster1971trivialen}. This is, after shrinking $\Delta$ if necessary we can assume that $\X\simeq X\times \Delta$ as deformations.

On the other hand, the rigidity of the singularities of $\F$ implies that $\dF$ is locally analytically trivial by the main theorem in \cite[page 31]{SuwaVersality}. In particular, making $\Delta$ smaller if necessary we can assume that there is a Stein cover $\U = \{U_i\}$ and vector fields $v_i\in T_\dF\vert_{U_i\times\Delta}$ projecting to $\frac{\partial}{\partial t}$ under the differential of $U_i\times\Delta \to \Delta$ 
such that integration along their flow gives the desired trivialization.
Observe that on each $U_{ij}$, the elements $v_{ij} = v_j-v_i$ project to zero and therefore these define a $\smash{\check{C}}$ech $1$-cocycle $\{v_{ij}\}\in H^1(\X,T_{\dF/\Delta})$.

The (local) triviality of $\F$ also shows that $\Sing(\dF)\to \Delta$ is locally trivial and hence flat. 
The exactness of the sequence 
$$ 0 \rightarrow T_{\dF/\Delta}\rightarrow T_{\X/\Delta}\rightarrow \I_{\Sing(\dF)}\otimes N_{\dF/\Delta}\rightarrow 0$$
now implies that $T_{\dF/\Delta}$ is also flat over $\Delta$. 
With this in mind, by the semicontinuity theorem and the hypothesis on $T_\F$ we can guarantee $H^1(\X,T_{\dF/\Delta})=0$. In particular, we can find vector fields $\overline{v}_i\in T_{\dF/\Delta}\vert_{U_i}$ such that  $v_{ij} = \overline{v}_j - \overline{v}_i$. The new family of vector fields $\{v_i-\overline{v}_i\}$ will now define a global section $v\in H^0(X,T_{\dF})$ which also projects to $\frac{\partial}{\partial t}$. To conclude, it suffices to observe that integrating the vector field $v$ with initial conditions $X\times\{0\}$ we obtain an automorphism of $X\times \Delta$ mapping the trivial unfolding to $\dF$.
\end{proof}

\bibliographystyle{alpha}
\bibliography{bibliografia.bib}

\end{document}